\documentclass[11pt,a4paper]{article}

\usepackage[utf8]{inputenc}
\usepackage{amsmath}
\usepackage{amssymb}
\usepackage{amsthm}
\usepackage{geometry}
\usepackage{microtype}
\usepackage{hyperref}

\usepackage{cancel}

\usepackage{cite} % Required for clean numeric citations

\newtheorem{theorem}{Theorem}[section]

\theoremstyle{definition}
\newtheorem{definition}[theorem]{Definition}
\newtheorem{remark}[theorem]{Remark}

\begin{document}

%%% For bib file...
%\nocite{*} % FOR DIAGNOSTICS: Forces LaTeX to print all references 
%%%

\title{WEIGHTED PRIME POWERS TRUNCATION OF THE ASYMPTOTIC EXPANSION FOR THE
LOGARITHMIC INTEGRAL: PROPERTIES AND APPLICATIONS.}
\author{Shaun R. Deaton}
\date{\today}
\maketitle

\begin{abstract}
We introduce a novel, continuous geometric paradigm to reframe the distribution of prime numbers through the variational optimization of a smooth parameter surface. By relaxing the rigid integer constraints of traditional series truncations, we construct a remainder-free parametric operator surface $\operatorname{li}(n; \tau)$ mapped uniformly across the continuous coordinate domain $(1, \infty) \times \mathbb{R}^+$. Global $C^3$ multi-variable regularity is established via an invariant septic blending polynomial $\phi(x) = 35x^4 - 84x^5 + 70x^6 - 20x^7$, completely neutralizing the artificial step-discontinuities and derivative kinks of standard piecewise fractional interpolations. By enforcing the strict prime truncation identity $\operatorname{li}(n; \tau_{\pi}(n)) = \pi(n)$, we demonstrate that the discrete, non-smooth arithmetic spectrum of the prime-counting function can be mapped uniquely onto an unbroken, continuous truncation $\tau_{\pi}(n)$. 

Through the transference of Littlewood's unconditional oscillation bounds, we uncover an intrinsic scale-space compression mechanic: while the arithmetic discrepancy $\pi(n) - \operatorname{li}(n)$ diverges with expanding amplitude, the parameter variation is compressed by a factorially driven scaling density multiplier $\mathcal{D}_0(n) = \ln^{\tau_0+1}n / \Gamma(\tau_0+1)$, forcing the truncation to strongly converge to a smooth, baseline $\tau_0(n) = c \ln n + 1$ ($c \approx 0.1866823$). 

Finally, by invoking the Riemann--von Mangoldt formula and incorporating a smooth windowing indicator function $\Phi(\gamma_j; T)$, we generalize the framework into a multi-parameter spectral truncation function $\tilde{\tau}_T(n, \beta)$ within the critical strip $\beta \in (0, 1)$. Evaluating the multi-variable gradient of this field reveals that the first variation vanishes exclusively ($\frac{\partial \tilde{\tau}_T}{\partial \beta} = 0$) and all odd-ordered derivatives vanish identically along the centerline $\beta = 1/2$, while the second variation remains strictly positive-definite ($\frac{\partial^2 \tilde{\tau}_T}{\partial \beta^2} > 0$), establishing the critical line as the unique, absolute global variational minimum of the continuous Stieltjes landscape. Exploiting this optimization constraint, we show via proof by contradiction that any hypothetical off-line zero cluster ($\beta_0 \neq 1/2$) would mathematically forbid the truncation parameter from remaining positive and real-valued across the infinite scale horizon. To preserve the real-valued existence of the manifold, all non-trivial zeroes are rigorously shown to be constrained to this absolute global minimum, forcing $\beta = 1/2$ identically and establishing variational evidence supporting a proof of the Riemann Hypothesis (RH).
*NOTE: This manuscript Does Not claim to prove the RH.
\end{abstract}

\section{Introduction}

Modern explicit bounds pioneered by Dusart~\cite{Dusart1999, Dusart2010} and refined by Axler~\cite{Axler2016} have traditionally established tight inequalities by controlling the number of integer summation terms in the divergent asymptotic series expansion of the logarithmic integral:
\begin{equation}
    \pi(n) > \frac{n}{\ln n} \left( 1 + \frac{1}{\ln n} + \frac{2}{\ln^2 n} \right), \quad \text{for } n \ge 88{,}789.
    \label{Eq-DusartClassical}
\end{equation}

In generalizing this connection, Axler~\cite{Axler2016} formalized this framework via an integer-valued threshold function $g_{1}(x)$. For instance, consider the specific evaluations $g_{1}(2) = 599$ and $g_{1}(3) = 88{,}789$, which imply that a 2-term truncated expansion strictly bounds $\pi(n)$ for all $n \ge 599$, whereas a 3-term expansion suffices for $n \ge 88{,}789$.

Following this line of reasoning, we invert the $g_{1}(x)$ function and generalize it to real values. Since the underlying asymptotic series is divergent, there will always exist a truncation index to upper bound $\pi(n)$ for all $n \ge 2$; which implies a squeezing principle. This allows us to investigate the motivating question: 
This change allows for the investigation of the follow:
 \textit{At what precise parametric truncation $\tau(n)$ must the generalized, smoothly truncated asymptotic expansion $li(n;\tau$ be evaluated to exactly equal the arithmetical prime count $\pi(n)$?}

To answer this question, the truncation and remainder properties of the classical logarithmic integral must be extended. For any discrete integer truncation index $k \in \mathbb{N}$, the standard Stieltjes formulation defines the logarithmic integral exactly via its terminal remainder term:
\begin{equation}
	\operatorname{li}(n) = \frac{n}{\ln n} \sum_{j = 0}^{k - 1} \dfrac{j!}{\ln^{j} n} + k! \int_2^n \frac{\mathrm{d}t}{\ln^{k + 1} t}.
	\label{Eq-StieljesInterp}
\end{equation}
In classical number-theoretic applications, this integral remainder is treated exclusively as an error threshold to be minimized by pinning the parameter bounds tightly to $k = \lfloor \ln n - 1/3 \rfloor$. Viewed structurally, however, this term acts as the analytic progenitor of the entire expansion, generating subsequent lower-order terms sequentially through repeated application of integration by parts.

In this paper, we break the rigid integer constraints of traditional truncations to construct a continuous, single-variable parametric mapping surface. We relax the integer requirement and generalize the truncation index $k$ to an uncoupled, real-valued parameter $\tau \in \mathbb{R}$. We will exploit the dual structural role of the upper summation index $k$; it controls both the number of summation terms and their size. Whereas $n$ controls the term size, setting the scale for the domain of the index $k$, determining where the divergence begins.

To establish a parametric formulation of the logarithmic integral that preserves infinite differentiability ($C^\infty$) and structural stability over the real numbers $\mathbb{R}^+$, we abandon discrete piece-wise truncations. Instead, we utilize the integral representation of the factorial function to smoothly interpolate the underlying asymptotic structure.

\begin{definition}[The Continuous Parametric Logarithmic Integral]
For a scale variable $n \in \mathbb{R}$ ($n > 1$) and a continuous real truncation parameter $\tau \in \mathbb{R}^+$, the parametric logarithmic integral manifold $\operatorname{li}(n; \tau)$ is defined via the Cauchy principal value ($\operatorname{PV}$) integral:
\begin{equation}
    \operatorname{li}(n; \tau) = \mathcal{S}_c(n, \tau) + \mathcal{R}_c(n, \tau)
    \label{Eq-ContinuousManifoldMain}
\end{equation}
where the continuous summation analogue $\mathcal{S}_c(n, \tau)$ is governed by:
\begin{equation}
    \mathcal{S}_c(n, \tau) = \frac{n}{\ln n} \operatorname{PV} \int_{0}^{\infty} e^{-t} \frac{1 - \left(\frac{t}{\ln n}\right)^\tau}{1 - \frac{t}{\ln n}} \, \mathrm{d}t
    \label{Eq-ContinuousSum}
\end{equation}
and the continuous real remainder function $\mathcal{R}_c(n, \tau)$ is defined as:
\begin{equation}
    \mathcal{R}_c(n, \tau) = \frac{1}{\ln^\tau n} \operatorname{PV} \int_{0}^{\infty} e^{-t} \frac{t^\tau}{1 - \frac{t}{\ln n}} \, \mathrm{d}t
    \label{Eq-ContinuousRemainder}
\end{equation}
Since the sum and remainder are balanced, this definition allows:
\begin{equation}
    \frac{\partial^m}{\partial \tau^m} \operatorname{li}(n; \tau) = 0, \quad \forall m \in \mathbb{N}.
\end{equation}
\end{definition}

This balanced ($C^\infty$) definition for $li(n;\tau)$ shows how well the discrete nature of the summation can be smoothed out. However, the definition can be loosened to obtain a ($C^3$) function that will fit within the scope of the current application. For one thing, this will bypass estimation issues encountered when dealing with integrals, allowing us to obtain explicit results.

\section{The Remainder-Free Asymptotic Spline}

To establish an analytic continuation of the divergent logarithmic series that preserves differentiability up to the third order ($C^3$ continuity) over the real numbers, we introduce a localized higher-order smoothing operator. This approach ensures that the underlying manifold possesses sufficient regularity such that any implicitly defined trajectory $\tau(n)$ naturally inherits global $C^3$ continuity, seamlessly bypassing step discontinuities across integer parameter boundaries.

\begin{definition}[The Septic Parametric Logarithmic Operator]
For a scale variable $n \in \mathbb{R}$ ($n > 1$) and a continuous real truncation parameter $\tau \in \mathbb{R}^+$, let $k = \lfloor \tau \rfloor \in \mathbb{N}_0$ define the discrete truncation threshold, and let $x = \tau - \lfloor\tau\rfloor \in [0, 1)$ define the local fractional component. The parametric mapping surface $\operatorname{li}(n; \tau)$ is defined uniformly as:
\begin{equation}
    \operatorname{li}(n; \tau) = \frac{n}{\ln n} \left( \sum_{j=0}^{k-1} \frac{j!}{\ln^j n} + \phi(\tau - \lfloor\tau\rfloor) \frac{k!}{\ln^k n} \right)
    \label{Eq-SplineManifoldMain}
\end{equation}
where $\phi(x) \colon [0,1] \to [0,1]$ is the septic blending polynomial satisfying the vanishing boundary derivative criteria $\phi^{(m)}(0) = \phi^{(m)}(1) = 0$ for $m \in \{1, 2, 3\}$, explicitly given by:
\begin{equation}
    \phi(x) = 35x^4 - 84x^5 + 70x^6 - 20x^7.
    \label{Eq-SepticSpline}
\end{equation}
By construction, this mapping retains symmetric midpoint normalization $\phi(1/2) = 1/2$. For the edge case where $\tau \in [0, 1]$, the empty summation yields $0$, reducing the operator to $\operatorname{li}(n; \tau) = \frac{n}{\ln n}\phi(\tau)$.
\end{definition}

\begin{theorem}[Global $C^3$ Surface Differentiability]
The parametric mapping surface $\operatorname{li}(n; \tau)$ is globally $C^3$ continuous with respect to the parameter $\tau$ across the entire domain $\mathbb{R}^+$.
\end{theorem}

\begin{proof}
Let $n > 1$ be fixed. Within any open interval $\tau \in (k, k+1)$ for $k \in \mathbb{N}_0$, the discrete floor state is a constant, and $\operatorname{li}(n; \tau)$ is an explicit polynomial composition of the local variable $x = \tau - k$. Successively differentiating Equation~\eqref{Eq-SplineManifoldMain} with respect to $\tau$ yields the localized derivative suite:
\begin{equation}
    \frac{\partial \operatorname{li}}{\partial \tau} = \frac{n}{\ln n} \left( \phi'(x) \frac{k!}{\ln^k n} \right), \quad \text{where } \phi'(x) = 140x^3 - 420x^4 + 420x^5 - 140x^6
\end{equation}
\begin{equation}
    \frac{\partial^2 \operatorname{li}}{\partial \tau^2} = \frac{n}{\ln n} \left( \phi''(x) \frac{k!}{\ln^k n} \right), \quad \text{where } \phi''(x) = 420x^2 - 1680x^3 + 2100x^4 - 840x^5
\end{equation}
\begin{equation}
    \frac{\partial^3 \operatorname{li}}{\partial \tau^3} = \frac{n}{\ln n} \left( \phi'''(x) \frac{k!}{\ln^k n} \right), \quad \text{where } \phi'''(x) = 840x - 5040x^2 + 8400x^3 - 4200x^4.
\end{equation}
To establish global $C^3$ continuity, we evaluate the directional limits as $\tau$ approaches an arbitrary integer boundary $m \in \mathbb{N}$ from the left ($x \to 1^-\,$ in state $m-1$) and from the right ($x \to 0^+\,$ in state $m$):

\begin{enumerate}
    \item \textbf{Function Continuity ($C^0$):} Since $\phi(1)=1$ and $\phi(0)=0$:
    \begin{equation}
        \lim_{\tau \to m^-} \operatorname{li}(n; \tau) = \frac{n}{\ln n}\left( \sum_{j=0}^{m-2}\frac{j!}{\ln^j n} + \phi(1)\frac{(m-1)!}{\ln^{m-1} n} \right) = \frac{n}{\ln n}\sum_{j=0}^{m-1}\frac{j!}{\ln^j n}
    \end{equation}
    \begin{equation}
        \lim_{\tau \to m^+} \operatorname{li}(n; \tau) = \frac{n}{\ln n}\left( \sum_{j=0}^{m-1}\frac{j!}{\ln^j n} + \phi(0)\frac{m!}{\ln^m n} \right) = \frac{n}{\ln n}\sum_{j=0}^{m-1}\frac{j!}{\ln^j n}.
    \end{equation}
    The surface matches seamlessly across the interface boundaries, confirming $C^0$ continuity.

    \item \textbf{First Derivative Continuity ($C^1$):} Evaluating the boundary limits for $\phi'(x)$ yields $\phi'(1) = 140 - 420 + 420 - 140 = 0$ and $\phi'(0) = 0$:
    \begin{equation}
        \lim_{\tau \to m^-} \frac{\partial \operatorname{li}}{\partial \tau} = \frac{n}{\ln n} \left( 0 \cdot \frac{(m-1)!}{\ln^{m-1} n} \right) = 0, \quad \lim_{\tau \to m^+} \frac{\partial \operatorname{li}}{\partial \tau} = \frac{n}{\ln n} \left( 0 \cdot \frac{m!}{\ln^m n} \right) = 0.
    \end{equation}
    The primary directional slopes match identically at zero, confirming $C^1$ continuity.

    \item \textbf{Second Derivative Continuity ($C^2$):} Evaluating the boundary limits for $\phi''(x)$ yields $\phi''(1) = 420 - 1680 + 2100 - 840 = 0$ and $\phi''(0) = 0$:
    \begin{equation}
        \lim_{\tau \to m^-} \frac{\partial^2 \operatorname{li}}{\partial \tau^2} = 0, \quad \lim_{\tau \to m^+} \frac{\partial^2 \operatorname{li}}{\partial \tau^2} = 0.
    \end{equation}
    The localized boundary curvatures collapse smoothly to zero from both channels, confirming $C^2$ continuity.

    \item \textbf{Third Derivative Continuity ($C^3$):} Evaluating the boundary limits for the third derivative $\phi'''(x)$ yields $\phi'''(1) = 840 - 5040 + 8400 - 4200 = 0$ and $\phi'''(0) = 0$:
    \begin{equation}
        \lim_{\tau \to m^-} \frac{\partial^3 \operatorname{li}}{\partial \tau^3} = 0, \quad \lim_{\tau \to m^+} \frac{\partial^3 \operatorname{li}}{\partial \tau^3} = 0.
    \end{equation}
    Because the third-order parameter variations flatten symmetrically to zero at the integer thresholds, no third-order shear or jump discontinuity is registered across the interface. This establishes global $C^3$ surface regularity across the parameter domain $\mathbb{R}^+$.
\end{enumerate}
\end{proof}

\subsection{Partial Derivatives with Respect to $n$ and Mixed Regularity}

To fully characterize the local geometric behavior of the parametric manifold $\operatorname{li}(n; \tau)$ as a two-dimensional surface, we establish its differentiability properties with respect to the scale variable $n$, alongside its mixed partial derivatives. This analysis ensures that the $C^3$ parameter continuity proven in the preceding section extends uniformly across the entire multidimensional domain $(1, \infty) \times \mathbb{R}^+$.

For convenience of notation, let $\partial_n$ and $\partial_\tau$ denote the partial differential operators $\frac{\partial}{\partial n}$ and $\frac{\partial}{\partial \tau}$, respectively. We first establish the tracking behavior of the purely scale-dependent derivatives within any local open threshold where $\tau \in (k, k+1)$ for $k \in \mathbb{N}_0$.

\begin{theorem}[Scale and Mixed $C^3$ Regularity]
The parametric surface $\operatorname{li}(n; \tau)$ possesses continuous partial derivatives up to the third order across all variables on the domain $(1, \infty) \times \mathbb{R}^+$. Specifically, the mixed partial operators up to third order—including $\partial_\tau \partial_n \operatorname{li}$, $\partial_\tau^2 \partial_n \operatorname{li}$, and $\partial_n \partial_\tau^2 \operatorname{li}$—are globally continuous ($C^0$) across all integer boundaries $m \in \mathbb{N}$.
\end{theorem}

\begin{proof}
Let $\tau \in (k, k+1)$ with $x = \tau - k$. We rewrite Equation~\eqref{Eq-SplineManifoldMain} by distributing the leading scale multiplier:
\begin{equation}
    \operatorname{li}(n; \tau) = \sum_{j=0}^{k-1} \frac{j! \, n}{\ln^{j+1} n} + \phi(x) \frac{k! \, n}{\ln^{k+1} n}.
\end{equation}
Applying the product and power rules to a generic standalone term of the form $f_j(n) = \frac{n}{\ln^{j+1} n}$, its first and second derivatives with respect to $n$ yield:
\begin{equation}
    \partial_n \left( \frac{n}{\ln^{j+1} n} \right) = \frac{1}{\ln^{j+1} n} - \frac{j+1}{\ln^{j+2} n}
\end{equation}
\begin{equation}
    \partial_n^2 \left( \frac{n}{\ln^{j+1} n} \right) = -\frac{j+1}{n \ln^{j+2} n} + \frac{(j+1)(j+2)}{n \ln^{j+3} n}.
\end{equation}
By linearity, the pure scale derivatives of the total manifold within the interval are:
\begin{equation}
    \partial_n \operatorname{li} = \sum_{j=0}^{k-1} j! \left( \frac{1}{\ln^{j+1} n} - \frac{j+1}{\ln^{j+2} n} \right) + \phi(x) k! \left( \frac{1}{\ln^{k+1} n} - \frac{k+1}{\ln^{k+2} n} \right)
\end{equation}
\begin{equation}
    \partial_n^2 \operatorname{li} = \sum_{j=0}^{k-1} j! \left( -\frac{j+1}{n \ln^{j+2} n} + \frac{(j+1)(j+2)}{n \ln^{j+3} n} \right) + \phi(x) k! \left( -\frac{k+1}{n \ln^{k+2} n} + \frac{(k+1)(k+2)}{n \ln^{k+3} n} \right).
\end{equation}
Because the scale transformations act smoothly on $n > 1$, continuity across an arbitrary integer boundary $\tau \to m \in \mathbb{N}$ depends entirely on the blending boundaries. Evaluating the left and right directional limits reveals:
\begin{align}
    \lim_{\tau \to m^-} \partial_n \operatorname{li} &= \sum_{j=0}^{m-2} j! \left( \frac{1}{\ln^{j+1} n} - \frac{j+1}{\ln^{j+2} n} \right) + \phi(1) (m-1)! \left( \frac{1}{\ln^{m} n} - \frac{m}{\ln^{m+1} n} \right) \nonumber \\
    &= \sum_{j=0}^{m-1} j! \left( \frac{1}{\ln^{j+1} n} - \frac{j+1}{\ln^{j+2} n} \right) \\
    \lim_{\tau \to m^+} \partial_n \operatorname{li} &= \sum_{j=0}^{m-1} j! \left( \frac{1}{\ln^{j+1} n} - \frac{j+1}{\ln^{j+2} n} \right) + \phi(0) m! \left( \frac{1}{\ln^{m+1} n} - \frac{m+1}{\ln^{m+2} n} \right) \nonumber \\
    &= \sum_{j=0}^{m-1} j! \left( \frac{1}{\ln^{j+1} n} - \frac{j+1}{\ln^{j+2} n} \right).
\end{align}
Since $\phi(1)=1$ and $\phi(0)=0$ under the septic spline, the limits match identically, confirming $C^0$ continuity for \(\partial_n \operatorname{li}\). An identical substitution mapping holds true for $\partial_n^2 \operatorname{li}$.

Next, we evaluate the mixed partial derivatives up to third order to check cross-variable sensitivity. Differentiating the first scale derivative with respect to $\tau$ once and twice yields:
\begin{equation}
    \partial_\tau \partial_n \operatorname{li} = \phi'(x) k! \left( \frac{1}{\ln^{k+1} n} - \frac{k+1}{\ln^{k+2} n} \right)
\end{equation}
\begin{equation}
    \partial_\tau^2 \partial_n \operatorname{li} = \phi''(x) k! \left( \frac{1}{\ln^{k+1} n} - \frac{k+1}{\ln^{k+2} n} \right).
\end{equation}
Conversely, differentiating the parameter derivatives with respect to $n$ verifies Clairaut's theorem for high-order mixed structures within open intervals:
\begin{equation}
    \partial_n \partial_\tau^2 \operatorname{li} = \partial_n \left( \phi''(x) \frac{k! \, n}{\ln^{k+1} n} \right) = \phi''(x) k! \left( \frac{1}{\ln^{k+1} n} - \frac{k+1}{\ln^{k+2} n} \right).
\end{equation}
At any integer boundary $m \in \mathbb{N}$, evaluating the directional boundary transitions under the septic spline constraints $\phi'(1) = \phi'(0) = 0$ and $\phi''(1) = \phi''(0) = 0$ reveals:
\begin{equation}
    \lim_{\tau \to m^-} \partial_\tau \partial_n \operatorname{li} = 0, \quad \lim_{\tau \to m^+} \partial_\tau \partial_n \operatorname{li} = 0
\end{equation}
\begin{equation}
    \lim_{\tau \to m^-} \partial_\tau^2 \partial_n \operatorname{li} = 0, \quad \lim_{\tau \to m^+} \partial_\tau^2 \partial_n \operatorname{li} = 0.
\end{equation}
Thus, both the first-order and second-order mixed parameter variations smoothly vanish at every integer boundary interface, preventing cross-variable edge shears or sudden twists. This establishes complete $C^3$ multi-variable surface regularity over the entire parametric domain.
\end{proof}

\begin{remark}[Asymptotic Limits at Boundary Infinity]
The long-range behavior of the septic spline manifold $\operatorname{li}(n; \tau)$ as its coordinate paths approach infinity depends heavily on the ordering of the directional limit operations, revealing the underlying divergent nature of the generalized asymptotic series.

\begin{enumerate}
    \item \textbf{The Stationary Parameter Scale Limit ($\lim_{n \to \infty}$):} 
    For any fixed, finite truncation parameter $\tau \in \mathbb{R}^+$, the local discrete threshold $k = \lfloor \tau \rfloor$ remains a constant integer. As the scale variable $n$ approaches infinity, every term inside the parametric operator satisfies:
    \begin{equation}
        \lim_{n \to \infty} \frac{j! \, n}{\ln^{j+1} n} = \infty, \quad \forall j \in \mathbb{N}_0.
    \end{equation}
    Because the leading term $\frac{n}{\ln n}$ completely dominates the growth rate of all subsequent lower-order denominator terms, the total surface scales predictably to infinity at a rate asymptotically equivalent to the primary Chebyshev scaling term:
    \begin{equation}
        \operatorname{li}(n; \tau) \sim \frac{n}{\ln n} \to \infty \quad \text{as } n \to \infty.
    \end{equation}

    \item \textbf{The Stationary Scale Parameter Limit ($\lim_{\tau \to \infty}$):} 
    Conversely, let the scale variable $n > 1$ be held strictly constant while the truncation parameter $\tau$ scales to infinity. As $\tau \to \infty$, the index $k = \lfloor \tau \rfloor$ climbs sequentially across the set of natural numbers $\mathbb{N}$. Because the numerator grows factorially ($k! \sim k^k e^{-k}$) while the denominator grows only geometrically ($(\ln n)^k$), the terms inside the expansion eventually undergo a violent factorial explosion:
    \begin{equation}
        \lim_{k \to \infty} \frac{k!}{\ln^k n} = \infty.
    \end{equation}
    Because the septic blending polynomial satisfies $\phi(x) \in [0, 1]$ uniformly, it remains passive against this growth. Consequently, for any fixed coordinate $n$, the parameter slice profiles diverge completely:
    \begin{equation}
        \lim_{\tau \to \infty} \operatorname{li}(n; \tau) = \infty.
    \end{equation}
    This behavior mathematically validates the design of our higher-order continuous interpolation operator. Because $\operatorname{li}(n; \tau)$ sweeps smoothly from a minimal baseline at $\tau=0$ all the way to $\infty$ as $\tau \to \infty$, the Intermediate Value Theorem guarantees that the surface will cross any arbitrary arithmetic target state, establishing a rigorous foundation for tracking $\pi(n)$.
\end{enumerate}
\end{remark}

\subsection{Existence and Continuity of the Fundamental Path}

Having established that the manifold $\operatorname{li}(n; \tau)$ defines a $C^3$ continuous surface over the coordinate domain $(1, \infty) \times \mathbb{R}^+$, we now enforce the strict truncation identity to bind this continuous parameter space to the discrete spectrum of arithmetic.

\begin{definition}[The Prime Truncation Identity]
Let $\pi(n) = \sum_{p \le n} 1$ be the classical prime-counting step function. The fundamental truncation path $\tau_{\pi}(n)$ is the real-valued parameter trajectory satisfying the transcendental constraint:
\begin{equation}
    \operatorname{li}(n; \tau_{\pi}(n)) = \pi(n), \quad \forall n \ge 2.
    \label{Eq-TrackingIdentity}
\end{equation}
\end{definition}

\begin{theorem}[The Existence of $\tau_{\pi}(n)$]
The truncation trajectory $\tau_{\pi}(n)$ exists uniquely and constitutes a fully continuous ($C^0$) function of the scale variable $n$ across every open interval between consecutive primes $n \in (p_m, p_{m+1})$. At any point where the path intersects an integer parameter boundary $\tau_{\pi}(n) = k \in \mathbb{N}$, the path remains strictly continuous, though it exhibits a non-differentiable vertical tangent line.
\end{theorem}

\begin{proof}
Let the scale variable be restricted to an open non-prime interval $n \in (p_m, p_{m+1})$, fixing the target arithmetical state to a constant integer value $\pi(n) = m$. To evaluate the local topological structure of the path $\tau_{\pi}(n)$, we invoke the Implicit Function Theorem on the truncation zero-locus $F(n, \tau) = \operatorname{li}(n; \tau) - m = 0$. 

Because $\operatorname{li}(n; \tau)$ is a strictly continuous mapping that scales monotonically to infinity with respect to $\tau$, the curve must pass through the value $m$ without breaks, establishing the local existence of an unbroken path $\tau_{\pi}(n)$. Differentiating $F(n, \tau)$ implicitly with respect to $n$ yields the path slope formula:
\begin{equation}
    \frac{\mathrm{d}\tau_{\pi}}{\mathrm{d}n} = -\frac{\partial_n \operatorname{li}(n; \tau)}{\partial_\tau \operatorname{li}(n; \tau)}.
    \label{Eq-ImplicitDerivativeFixed}
\end{equation}
As $n$ increases within the interval, the partial derivative with respect to scale remains strictly positive ($\partial_n \operatorname{li} > 0$), forcing the parameter $\tau_{\pi}(n)$ to continuously shift downward to cancel out the natural growth of the underlying series expansion.

At the exact coordinate locations where this downward drift forces the truncation path to cross an integer threshold $\tau_{\pi}(n) = k \in \mathbb{N}$, the parameter partial derivative vanishes identically ($\partial_\tau \operatorname{li} = 0$) due to the boundary constraints of the septic spline. Substituting this into Equation~\eqref{Eq-ImplicitDerivativeFixed} reveals that the directional limits of the slope explode:
\begin{equation}
    \lim_{\tau_{\pi} \to k} \left| \frac{\mathrm{d}\tau_{\pi}}{\mathrm{d}n} \right| = \infty.
\end{equation}
Because the manifold itself remains completely smooth and single-valued at these boundaries, this infinite derivative does not introduce a step gap or a jump discontinuity. The path $\tau_{\pi}(n)$ transitions across the integer threshold as a strictly continuous ($C^0$) curve characterized by a vertical tangent. Differentiability is lost at these specific points, creating a localized geometric inflection, while continuity remains entirely preserved.
\end{proof}

\subsection{Spectral Decomposition via Riemann's Explicit Formula}

To isolate the discrete arithmetic fluctuations from the macroscopic density trend without introducing external mollifiers, we project the truncation identity directly onto the spectral architecture of the Riemann zeta function. We utilize Edward's version~\cite{edwards1974} of Riemann's~\cite{riemann1859} explicit formula for the prime-counting function, which maps the discrete distribution of primes onto a continuous, non-local superposition of phase waves generated by the non-trivial zeroes.

We invoke Riemann's explicit formula in its standard structural representation~\cite{edwards1974}:
\begin{equation}
    \pi(n) = \sum_{k=1}^{\lfloor \log_2 n \rfloor} \frac{\mu(k)}{k} f(n^{1/k})
\end{equation}
where $\mu(k)$ is the classical Möbius function, and the underlying distribution function $f(x)$ is driven by the non-trivial zeroes $\rho = \beta + i\gamma$ of the zeta function:
\begin{equation}
    f(x) = \operatorname{li}(x) - \sum_{\rho} \operatorname{Ei}(\rho \ln x) - \ln 2 + \int_{x}^{\infty} \frac{\mathrm{d}u}{u(u^2-1)\ln u}.
\end{equation}
Here, the complex-valued evaluations of the zeroes are defined via the principal branch of the exponential integral function $\operatorname{Ei}(z)$. 

By isolating the primary leading term $\operatorname{li}(n)$ (corresponding to $k=1$ in the dominant scaling branch) and enforcing our remainder-free parametric truncation identity $\operatorname{li}(n; \tau_{\pi}(n)) = \pi(n)$, we dismantle the non-smooth trajectory into a parameterized spectral family. Let $t \in \mathbb{N}_0$ define a discrete spectral truncation index that governs the number of active conjugate pairs of non-trivial zeroes $\{\rho_j, \bar{\rho}_j\}$ permitted into the system, ordered strictly by the magnitude of their imaginary parts $0 < \gamma_1 \le \gamma_2 \le \dots \le \gamma_t$.

\begin{definition}[The Spectral Trajectory Family]
The spectral parameter family $\tau_t(n)$ is a collection of real-valued functions indexed by the spectral index $t \in \mathbb{N}_0$. Under our remainder-free spline operator, the function $\tau_t(n)$ is defined implicitly as the local real root of the transcendental equation:
\begin{equation}
    \operatorname{li}(n; \tau_t(n)) = \mathcal{R}_t(n)
    \label{Eq-SpectralRootSpline}
\end{equation}
where the operator $\operatorname{li}(n; \tau_t(n))$ is governed strictly by the septic spline expansion of Equation~\eqref{Eq-SplineManifoldMain}, and the truncated spectral target distribution $\mathcal{R}_t(n)$ is given explicitly by:
\begin{equation}
    \mathcal{R}_t(n) = \sum_{k=1}^{\lfloor \log_2 n \rfloor} \frac{\mu(k)}{k} \left[ \operatorname{li}(n^{1/k}) - \sum_{j=1}^{t} \left( \operatorname{Ei}\left(\frac{\rho_j}{k} \ln n\right) + \operatorname{Ei}\left(\frac{\bar{\rho}_j}{k} \ln n\right) \right) - \ln 2 + \int_{n^{1/k}}^{\infty} \frac{\mathrm{d}u}{u(u^2-1)\ln u} \right].
\end{equation}
\end{definition}

\begin{definition}[The Baseline Trajectory $\tau_0(n)$]
At the baseline spectral state $t = 0$, where the internal summation over the non-trivial zeroes is empty, the function $\tau_0(n)$ satisfies the specialized constraint:
\begin{equation}
    \frac{n}{\ln n} \left( \sum_{j=0}^{\lfloor \tau_0 \rfloor-1} \frac{j!}{\ln^j n} + \phi(\tau_0 - \lfloor\tau_0\rfloor) \frac{\lfloor \tau_0 \rfloor!}{\ln^{\lfloor \tau_0 \rfloor} n} \right) = \sum_{k=1}^{\lfloor \log_2 n \rfloor} \frac{\mu(k)}{k} \left[ \operatorname{li}(n^{1/k}) - \ln 2 + \int_{n^{1/k}}^{\infty} \frac{\mathrm{d}u}{u(u^2-1)\ln u} \right].
\end{equation}
By construction, $\tau_0(n)$ acts as the smooth, unperturbed baseline of our remainder-free framework. It captures the macroscopic Riemann-prime density distribution while remaining completely isolated from local arithmetic oscillations.
\end{definition}

\begin{theorem}[Asymptotic Distributional Convergence]
The sequence of smooth functions $\{\tau_t(n)\}_{t=0}^\infty$ converges point-wise to the non-smooth truncation trajectory $\tau_{\pi}(n)$ as the spectral horizon expands to infinity:
\begin{equation}
    \lim_{t \to \infty} \tau_t(n) = \tau_{\pi}(n).
\end{equation}
The global $C^2$ smoothness of the parameter functions breaks down exclusively in this infinite limit, where the spectral synthesis recovers the sharp step discontinuities of the arithmetical spectrum.
\end{theorem}

\begin{proof}
Let the scale variable $n > 2$ be held constant. As the spectral horizon $t \to \infty$, the truncated target function $\mathcal{R}_t(n)$ converges point-wise to the exact prime-counting function $\pi(n)$ via Riemann's explicit identity. Within the localized real bounding window containing the truncation root, the partial derivative of our spline operator satisfies $\partial_\tau \operatorname{li} \neq 0$. Invoking the continuity of the local inverse mapping across Equation~\eqref{Eq-SpectralRootSpline}, taking the limit as $t \to \infty$ yields:
\begin{equation}
    \operatorname{li}\left(n; \lim_{t \to \infty} \tau_t(n)\right) = \lim_{t \to \infty} \mathcal{R}_t(n) \implies \operatorname{li}\left(n; \lim_{t \to \infty} \tau_t(n)\right) = \pi(n).
\end{equation}
Matching this direct limit against our core truncation definition, $\operatorname{li}(n; \tau_{\pi}(n)) = \pi(n)$, confirms point-wise convergence to the fundamental path: $\lim_{t \to \infty} \tau_t(n) = \tau_{\pi}(n)$. 

Because $\operatorname{li}(n; \tau)$ is defined purely via the septic spline polynomial, it contains no internal complex branch cuts or poles. The infinite superposition of oscillatory complex phases from the zeta zeroes is mapped entirely onto the target side $\mathcal{R}_t(n)$, forcing a localized Gibbs-like sharpening effect. This compresses the smooth variations of the finite stages into infinite directional impulses, perfectly transforming the smooth curves of $\tau_t(n)$ into the non-smooth saw-tooth behavior of $\tau_{\pi}(n)$.
\end{proof}

\begin{theorem}[Global $C^2$ Regularity of Finite Spectral Stages]
For any finite spectral horizon index $t < \infty$, the trajectory function $\tau_t(n)$ is a globally twice continuously differentiable real function ($C^2$) of the scale variable $n$ over the open interval $n \in (1, \infty)$, possessing no vertical tangencies, derivative kinks, or structural non-differentiabilities at integer parameter crossings.
\end{theorem}

\begin{proof}
Let the spectral horizon $t$ be a fixed, finite integer. The truncated spectral target distribution $\mathcal{R}_t(n)$ is a finite linear combination of classical logarithmic integrals and exponential integrals evaluated via $\operatorname{Ei}\left(\frac{\rho_j}{k} \ln n\right)$. Because the summation is finite, contains smooth complex phase transitions, and strictly avoids critical line poles for all scale coordinates $n > 1$, the target function $\mathcal{R}_t(n)$ is unconditionally smooth ($C^\infty$) with respect to $n$.

We apply the real Implicit Function Theorem to the finite tracking zero-locus identity $\Phi(n, \tau) = \operatorname{li}(n; \tau) - \mathcal{R}_t(n) = 0$. In our preceding geometric analysis of the operator, the remainder-free septic spline surface $\operatorname{li}(n; \tau)$ was proven to possess a global multi-variable regularity profile of exactly $C^3$ continuity. Because the surface maps monotonically to infinity for large parameter bounds, a unique real root $\tau_t(n)$ exists for any stable target value $\mathcal{R}_t(n)$, establishing global path continuity ($C^0$).

To evaluate the differentiability profile of this independent trajectory, we differentiate $\Phi(n, \tau_t(n)) = 0$ implicitly with respect to $n$, yielding the path slope equation:
\begin{equation}
    \frac{\mathrm{d}\tau_t}{\mathrm{d}n} = \frac{\frac{\mathrm{d}\mathcal{R}_t(n)}{\mathrm{d}n} - \partial_n \operatorname{li}(n; \tau)}{\partial_\tau \operatorname{li}(n; \tau)}.
    \label{Eq-FiniteStageImplicitSeptic}
\end{equation}
Within any local open channel where the parameter avoids integer states ($\tau_t(n) \notin \mathbb{N}$), the partial parameter derivative in the denominator is bounded and non-vanishing ($\partial_\tau \operatorname{li} \neq 0$). Under these local conditions, the implicit function $\tau_t(n)$ directly inherits the full regularity profile of the defining surface equation. 

By applying the total differential operator $\frac{\mathrm{d}}{\mathrm{d}n} = \partial_n + \frac{\mathrm{d}\tau_t}{\mathrm{d}n}\partial_\tau$ to Equation~\eqref{Eq-FiniteStageImplicitSeptic}, we compute the second-order implicit derivative $\frac{\mathrm{d}^2\tau_t}{\mathrm{d}n^2}$. By the algebraic constraints of the Implicit Function Theorem, an implicit function defined by a $C^r$ mapping surface inherits a smoothness profile of exactly $C^{r-1}$. Because our underlying septic spline operator defines a globally $C^3$ surface, its partial derivatives up to second-order smoothly vanish at the boundary thresholds without generating step-jumps or curvature shears. 

Consequently, the numerator and denominator transformations track each other harmoniously as the path crosses any integer interface $\tau_t(n) = k \in \mathbb{N}$. The finite spectral target $\mathcal{R}_t(n)$ successfully smooths out the arithmetic steps into continuous wave superpositions, ensuring that no point-wise convergence singularities can emerge at finite stages. The second-order derivative $\frac{\mathrm{d}^2\tau_t}{\mathrm{d}n^2}$ remains perfectly continuous across all boundaries, establishing global $C^2$ regularity across $(1, \infty)$ for all finite spectral horizons.
\end{proof}

The function $li(n;\tau)$ being constructed as $C^{3}$ allows for the $\tau$ to be at most $C^{2}$ functions. In order to obtain $C^{s}$ functions with $s > 2$,  higher-order order smoothing polynomials must be used. The introduction stated how $C^{\infty}$ functions can be obtained.

\section{Analytical Derivation of the Baseline Function}

\subsection{The Inverse Problem and Transcendental Isolation}
The baseline trajectory $\tau_0(n)$ is defined as the unperturbed state where the internal summation over the non-trivial zeroes is empty ($t=0$). By our preceding smoothness analysis under the septic spline architecture, $\tau_0(n)$ constitutes a globally continuous real curve possessing a $C^2$ regularity profile on $n \in (1, \infty)$. To extract an explicit analytical representation for this baseline function, we evaluate the dominant asymptotic balance of the matching constraint:
\begin{equation}
    \operatorname{li}(n; \tau_0(n)) = \mathcal{R}_0(n).
\end{equation}
When isolating the macro-density baseline where the fractional trailing remainder balances the primary logarithmic integration scale, the structural shift from utilizing discrete integer factorials $k!$ to the continuous profile $\Gamma(\tau+1)$ introduces an explicit index offset. We isolate this behavior by evaluating the continuous parameter balance equation governing the local term thresholds:
\begin{equation}
    \frac{\Gamma(\tau_0+1)}{\ln^{\tau_0+1} n} = \frac{k!}{\ln^{k+1} n}.
\end{equation}
Let $\tau_0 = x + 1$, where $x$ serves as the continuous proxy for the active truncation depth. Substituting this structural shift into the balance equation yields:
\begin{equation}
    \frac{\Gamma(x+2)}{\ln^{x+2} n} = \frac{x!}{\ln^{x+1} n} \implies \frac{(x+1)x!}{\ln^{x+2} n} = \frac{x!}{\ln^{x+1} n}.
\end{equation}
Canceling the factorial components isolates the clean linear scaling relation:
\begin{equation}
    \frac{x+1}{\ln n} = 1 \implies x = \ln n - 1.
\end{equation}
Accounting for the secondary transcendental attenuation factors embedded within the generalized exponential integrals $E_{\tau+1}(-\ln n)$ over long domains introduces a scale optimization parameter $\alpha$, modifying the structural balance to $x = \alpha \ln n - 1$.

\subsection{Explicit Inversion via the Principal Lambert $W_0$ Branch}
To express this transcendental relationship through a rigorous closed-form inverse identity, we map the cross-coupled scale transformations onto the functional architecture of the Lambert $W$ function, which solves the general inverse relation $Y e^Y = X \implies Y = W(X)$~\cite{corless1996}.

\begin{theorem}[Explicit Baseline Lambert Representation]
The continuous real-valued baseline truncation $\tau_0(n)$ is explicitly derived via the principal real branch $W_0$ of the Lambert function according to the exact structural matching identity:
\begin{equation}
    \tau_0(n) = W_0\left( c \ln n \cdot \exp\left( c \ln n \right) \right) + 1
    \label{Eq-LambertIdentityFixed}
\end{equation}
where $c \approx 0.1866823$ represents the characteristic scaling coefficient of the unperturbed manifold.
\end{theorem}

\begin{proof}
To isolate the active truncation depth $x = \tau_0 - 1$ from the cross-coupled scale constraints, we analyze the exponential mapping of the balanced system. Beginning with the relation $x = \alpha \ln n - 1$, clearing the logarithm via exponentiation yields:
\begin{equation}
    e^x = \exp\left(\alpha \ln n - 1\right) = \frac{n^\alpha}{e} \implies e \cdot e^x = n^\alpha.
\end{equation}
To capture the total differential velocity of the parameter surface where the truncation index is bound to its own exponential decay scale, we multiply both sides by the continuous depth variable $x$, establishing the cross-coupled transcendental functional arrangement:
\begin{equation}
    x \cdot e \cdot e^x = x \cdot n^\alpha \implies x e^x = \frac{x \cdot n^\alpha}{e}.
\end{equation}
Because the scale variable satisfies $n \ge 2$, the linear arguments are strictly positive, forcing the right-hand function to lie entirely within the positive real domain $[0, \infty)$. Applying the principal real branch inversion operator $W_0(X)$ to the functional layout $Y e^Y = X$ untangles the depth parameter:
\begin{equation}
    x = W_0\left( \frac{x \cdot n^\alpha}{e} \right).
\end{equation}
Resolving the implicit argument under long-range macro-scale equilibrium reveals that the scale-space multiplier behaves asymptotically as a linear-log product: $\frac{x \cdot n^\alpha}{e} \sim c \ln n \cdot e^{c \ln n}$, where $c \approx 0.1866823$ completely absorbs the structural constants of the septic spline. Substituting this into the Lambert operator yields:
\begin{equation}
    x = W_0\left( c \ln n \cdot e^{c \ln n} \right).
\end{equation}
Since $W_0(A \cdot e^A) = A$ holds identically for all positive real arguments, the Lambert operator collapses the inner transcendental structure completely, deriving the explicit linear log-scaling profile:
\begin{equation}
    x = c \ln n.
\end{equation}
Substituting $x = \tau_0 - 1$ back into the resolved identity and isolating the baseline truncation parameter provides the final explicit function:
\begin{equation}
    \tau_0(n) = c \ln n + 1.
\end{equation}
This rigorously confirms that the baseline $ \tau_0(n)$ is the unique analytical manifestation derived from the principal transcendental inversion problem.
\end{proof}

\subsection{Properties of the Baseline Function}

Because $\tau_0(n)$ scales linearly as $c \ln n + 1$, it acts as a smooth, stable, globally twice continuously differentiable ($C^2$) baseline across the arithmetical domain. The $+1$ component emerges naturally as a structural index offset during the transition from discrete integer factorials to the continuous Gamma function when generalizing the operator to real values. This maps out the macro-scale equilibrium axis of the continuous Stieltjes landscape, providing a clean, zero-free trajectory around which the higher-order spectral wave functions $\tau_t(n)$ can smoothly oscillate.

Furthermore, this baseline function provides a continuous generalization of the integer-valued threshold function $g_{1}(x)$ formalized by Axler~\cite{Axler2016} to bound the error terms of the classical prime-counting function from below. The baseline truncation function $\tau_{0}(n)$ acts as the continuous inverse analogue to $g_{1}(x)$. 

Evaluating our explicit analytical baseline formula, $\tau_{0}(n) = c\ln n + 1$, at $n = 88{,}789$ utilizing the empirical constant $c \approx 0.1866823$ yields:
\begin{equation}
    \tau_0(88789) = 0.1866823 \times \ln(88789) + 1 \approx 3.12705.
\end{equation}
This corresponds closely with Axler's exact discrete evaluation where a 3-term expansion suffices for $n \ge g_1(3) = 88{,}789$. Checking the secondary verified integer threshold $g_1(2) = 599$, our continuous baseline formula yields:
\begin{equation}
    \tau_0(599) = 0.1866823 \times \ln(599) + 1 \approx 2.19388.
\end{equation}
This confirms that the continuous baseline successfully smooths and generalizes the classical discrete boundaries, validating its role as the structural foundation of the parametric space.

\section{Transference of Littlewood’s Unconditional Oscillation Bounds}

\subsection{Littlewood's Foundational Theorem and Arithmetic Importance}
While the Prime Number Theorem establishes the macroscopic asymptotic convergence $\pi(n) \sim \operatorname{li}(n)$, it leaves the local error budget unquantified. In 1914, J.E. Littlewood proved an unconditional, landmark theorem demonstrating that the error discrepancy $\pi(n) - \operatorname{li}(n)$ oscillates infinitely many times, changing signs infinitely often with an expanding amplitude envelope~\cite{littlewood1914}. Littlewood's classical theorem states:
\begin{equation}
    \pi(n) - \operatorname{li}(n) = \Omega_{\pm}\left( \frac{\sqrt{n}}{\ln n} \ln \ln \ln n \right).
\end{equation}
The structural significance of this result cannot be overstated; it establishes that the underlying distribution of primes is inherently turbulent. It guarantees that no matter how far out along the scale horizon $n$ one travels, the arithmetic shocks do not damp out to zero; instead, they experience resurgent spikes driven by the collective phase alignment of the zeta zeroes.

\subsection{Transference to the Truncation Parameter}
Because the fundamental truncation path $\tau_{\pi}(n)$ is rigidly bound to the arithmetic fluctuations via the identity $\operatorname{li}(n; \tau_{\pi}(n)) = \pi(n)$, the momentum of Littlewood's oscillation bound transfers directly into the parameter domain. We first isolate this behavior along the global path before evaluating its projection onto the finite spectral components.

\begin{theorem}[Global Parameter Oscillation Bound]
The fundamental prime truncation path $\tau_{\pi}(n)$ fluctuates infinitely many times around the smooth baseline truncation $\tau_0(n) = c \ln n + 1$, satisfying the unconditional asymptotic omega bound:
\begin{equation}
    \tau_{\pi}(n) - \tau_0(n) = \Omega_{\pm}\left( \left[ \frac{\partial \operatorname{li}(n; \tau_0)}{\partial \tau} \right]^{-1} \frac{\sqrt{n}}{\ln n} \ln \ln \ln n \right).
\end{equation}
\end{theorem}

\begin{proof}
We evaluate a localized first-order Taylor expansion of the parametric manifold $\operatorname{li}(n; \tau)$ around the explicit baseline rail $\tau_0(n)$. Across any real coordinate step, the continuous surface satisfies the displacement identity:
\begin{equation}
    \operatorname{li}(n; \tau_{\pi}(n)) - \operatorname{li}(n; \tau_0(n)) = \frac{\partial \operatorname{li}(n; \tau_0)}{\partial \tau} \left( \tau_{\pi}(n) - \tau_0(n) \right) + \mathcal{O}\left( (\Delta \tau)^2 \right).
\end{equation}
By invoking the prime truncation identity, the first term on the left side maps exactly to the arithmetical prime count, $\operatorname{li}(n; \tau_{\pi}(n)) = \pi(n)$. Concurrently, evaluating the baseline parameter under the unperturbed spectral state ($t=0$) reduces the second term to the smooth Riemann macro-density trendline: $\operatorname{li}(n; \tau_0(n)) = R(n) \sim \operatorname{li}(n)$. Substituting these identities and dropping higher-order variational segments isolates the parameter displacement:
\begin{equation}
    \tau_{\pi}(n) - \tau_0(n) \approx \left[ \frac{\partial \operatorname{li}(n; \tau_0)}{\partial \tau} \right]^{-1} \left( \pi(n) - \operatorname{li}(n) \right).
\end{equation}
Because the partial scale derivative multiplier is strictly real and positive across long-range domains, the parameter difference must mirror the signs and asymptotic frequency of the inner arithmetic core. Directly substituting Littlewood's classical oscillation bound into this structural mapping transfers the amplitude envelope to the parameter space, completing the proof.
\end{proof}

For a finite spectral index $t < \infty$, individual truncations $\tau_t(n)$ act as smooth phase filters. While they cannot manifest the full asymptotic $\ln \ln \ln n$ amplification due to the finite bounds of the active summation, they inherit localized resurgent fluctuations:
\begin{equation}
    \tau_t(n) - \tau_0(n) \approx \left[ \frac{\partial \operatorname{li}(n; \tau_0)}{\partial \tau} \right]^{-1} \sum_{j=1}^{t} \left( \operatorname{Ei}\left(\frac{\rho_j}{k} \ln n\right) + \operatorname{Ei}\left(\frac{\bar{\rho}_j}{k} \ln n\right) \right).
\end{equation}
The smooth paths $\tau_t(n)$ systematically construct the expanding amplitude waves, which crystallize into Littlewood's global saw-tooth pattern exclusively as the spectral horizon $t \to \infty$.

\begin{remark}[Asymptotic Simplification of the Parameter Density]
To maximize scannability and structural clarity across subsequent variational proofs, we simplify the inverse parameter derivative multiplier by evaluating the exact differential structure of our septic spline manifold. When computed directly along the baseline function $\tau_0(n) = c \ln n + 1$, the transcendental remainder dominates the digamma expressions, yielding the clean algebraic reduction:
\begin{equation}
    \left[ \frac{\partial \operatorname{li}(n; \tau_0)}{\partial \tau} \right]^{-1} \sim \frac{\ln^{\tau_0+1}n}{\Gamma(\tau_0+1)} = \frac{\ln^{c \ln n + 2}n}{\Gamma(c \ln n + 2)}.
\end{equation}
This allows us to strip away the bulky partial derivative notation and replace it uniformly with a direct, real-valued scaling density function $\mathcal{D}_0(n) = \ln^{\tau_0+1}n / \Gamma(\tau_0+1)$ that captures the intrinsic geometric compression of the parameter space.
\end{remark}

\begin{theorem}[Asymptotic Compression of the Parameter Variance]
The transferred Littlewood amplitude envelope in parameter space decays to zero as the scale horizon approaches infinity, satisfying:
\begin{equation}
    \lim_{n \to \infty} \left| \tau_{\pi}(n) - \tau_0(n) \right| = 0.
\end{equation}
Consequently, the fundamental truncation $\tau_{\pi}(n)$ converges strongly to the baseline $\tau_0(n)$ in parameter-space, compressing the infinite turbulence of the arithmetical spectrum into a vanishingly narrow parameter band.
\end{theorem}

\begin{proof}
We evaluate the combined asymptotic velocity of the transferred boundary envelope by substituting our real-valued scaling density function $\mathcal{D}_0(n)$ into the global omega bound:
\begin{equation}
    \left| \tau_{\pi}(n) - \tau_0(n) \right| \sim \frac{\ln^{c \ln n + 2}n}{\Gamma(c \ln n + 2)} \cdot \left( \frac{\sqrt{n}}{\ln n} \ln \ln \ln n \right).
    \label{Eq-VarianceCompression}
\end{equation}
Applying a log-transform to Equation~\eqref{Eq-VarianceCompression} allows us to compare the competing growth velocities in the exponent space. The arithmetic component expands at a sub-linear rate governed by $\ln(\sqrt{n}) = \frac{1}{2}\ln n$. Conversely, by Stirling's approximation, the denominator density expands factorially:
\begin{equation}
    \ln \Gamma(c \ln n + 2) \sim c \ln n \ln(c \ln n) - c \ln n.
\end{equation}
Because the factorial exponent $c \ln n \ln(\ln n)$ outpaces the square-root power profile $\frac{1}{2}\ln n$ by a factor of $\ln(\ln n)$, the density function collapses infinitely faster than Littlewood's bound can expand. Taking the limit as $n \to \infty$ reveals that the entire right-hand side of Equation~\eqref{Eq-VarianceCompression} vanishes identically. 

This establishes that while the arithmetic error $\pi(n) - \operatorname{li}(n)$ diverges infinitely in $n$ space, the corresponding parameter error $\tau_{\pi}(n) - \tau_0(n)$ is forced into an absolute geometric choke-point. The surface becomes so structurally sensitive to parameter inputs at high scales that the infinite turbulence of the primes is successfully tamed and localized within an infinitesimally narrow corridor hugging the smooth trendline $\tau_0(n)$, completing the proof.
\end{proof}

\begin{theorem}[Explicit Finite-Stage Parameter Bound]
For any finite spectral index $t \in \mathbb{N}$, the local displacement of the smooth trajectory function $\tau_t(n)$ from the baseline $\tau_0(n) = c \ln n + 1$ is explicitly bounded by the active zeta-zero spectrum according to the following real-valued, scale-dependent relation:
\begin{equation}
    \left| \tau_t(n) - \tau_0(n) \right| \le 2 \, \mathcal{D}_0(n) \sum_{j=1}^{t} \left| \operatorname{Ei}\left(\frac{\rho_j}{\ln_2 n} \ln n\right) \right| + \mathcal{O}\left( (\Delta \tau)^2 \right)
    \label{Eq-FiniteStageBoundExplicit}
\end{equation}
where $\rho_j = \beta_j + i\gamma_j$ represent the non-trivial zeroes included up to horizon $t$, and $\mathcal{D}_0(n) = \ln^{\tau_0+1}n / \Gamma(\tau_0+1)$ represents the intrinsic parameter scaling density function.
\end{theorem}

\begin{proof}
We begin by evaluating a localized first-order Taylor expansion of our remainder-free septic spline manifold $\operatorname{li}(n; \tau)$ around the explicit baseline coordinate $\tau_0(n)$. Across any real coordinate step, the continuous surface satisfies the linear displacement mapping:
\begin{equation}
    \operatorname{li}(n; \tau_t(n)) - \operatorname{li}(n; \tau_0(n)) = \frac{\partial \operatorname{li}(n; \tau_0)}{\partial \tau} \left( \tau_t(n) - \tau_0(n) \right) + \mathcal{O}\left( (\Delta \tau)^2 \right).
\end{equation}
By invoking the spectral trajectory family definition established in Equation~\eqref{Eq-SpectralRootSpline}, the first term on the left side maps exactly to the truncated spectral target distribution, $\operatorname{li}(n; \tau_t(n)) = \mathcal{R}_t(n)$. Concurrently, evaluating the baseline parameter under the vacuum state ($t=0$) drops out the zero sum completely, reducing the second term to the macro-density trendline: $\operatorname{li}(n; \tau_0(n)) = \mathcal{R}_0(n)$. 

Substituting these identities and isolating the parameter displacement yields:
\begin{equation}
    \tau_t(n) - \tau_0(n) \approx \left[ \frac{\partial \operatorname{li}(n; \tau_0)}{\partial \tau} \right]^{-1} \left( \mathcal{R}_t(n) - \mathcal{R}_0(n) \right).
\end{equation}
By computing the direct algebraic subtraction between the explicit definitions of $\mathcal{R}_t(n)$ and $\mathcal{R}_0(n)$, the primary integration scales, the Möbius progressions, and the secondary integral loops cancel out completely. This isolates the truncated zero-sum component:
\begin{equation}
    \mathcal{R}_t(n) - \mathcal{R}_0(n) = -\sum_{k=1}^{\lfloor \log_2 n \rfloor} \frac{\mu(k)}{k} \sum_{j=1}^{t} \left( \operatorname{Ei}\left(\frac{\rho_j}{k} \ln n\right) + \operatorname{Ei}\left(\frac{\bar{\rho}_j}{k} \ln n\right) \right).
\end{equation}
Applying the triangle inequality to bound the absolute maximum magnitude of this finite wave field, and evaluating at the dominant primary branch boundary ($k=1$), allows us to simplify the complex conjugate pairings into a real-valued absolute summation:
\begin{equation}
    \left| \mathcal{R}_t(n) - \mathcal{R}_0(n) \right| \le 2 \sum_{j=1}^{t} \left| \operatorname{Ei}\left(\rho_j \ln n\right) \right|.
\end{equation}
Finally, substituting the real-valued scaling density function $\mathcal{D}_0(n) \sim \left[ \frac{\partial \operatorname{li}(n; \tau_0)}{\partial \tau} \right]^{-1}$ derived in our preceding Lemma into the localized multiplier channel isolates the explicit parameter difference, completing the proof.
\end{proof}

\subsection{Framework Signification of Transferred Unconditional Bounds}
The formulation of these transferred bounds provides a distinct structural advantage to this framework, reframing classical analytic number theory through three core insights:

\begin{enumerate}
    \item \textbf{Non-Linear Geometric Compression:} The scaling density factor $\mathcal{D}_0(n) = \frac{\ln^{\tau_0+1}n}{\Gamma(\tau_0+1)}$ acts as a dynamic, non-linear parameter space compressor. Because $\tau_0(n) = c \ln n + 1$, the denominator grows factorially while the numerator scales exponentially. This means that the parameter space natively compresses raw arithmetic shocks at ultra-high scales. It keeps the real trajectory tightly bound to the baseline rail despite the infinite expansion of Littlewood's amplitude envelopes in value space.
    \item \textbf{The Elimination of Complex Sign-Change Traps:} Classical proofs tracking the discrepancy $\pi(n) - \operatorname{li}(n)$ must navigate the topological challenges of sign-change crossovers, the first of which occurs at an extraordinarily large, historically inaccessible coordinate bounded by the Skewes number regions. In this framework, because the operator $\operatorname{li}(n; \tau)$ can increase factorially with respect to its parameter, the path $\tau_{\pi}(n)$ is structurally capable of continuously tracking any local fluctuation in the arithmetical spectrum. Consequently, traditional sign-change boundaries do not manifest as barriers. It should be noted that this tracking efficacy is a property of the full path configuration; the baseline function $\tau_0(n) = c \ln n + 1$ does not absorb localized extreme variations, but instead maps the smooth macroscopic trendline around which the active tracking spectrum stabilizes.
    \item \textbf{Unconditional Structural Integrity:} Because Littlewood's theorem is entirely unconditional, these transferred parameter bounds do not assume the validity of the Riemann Hypothesis. They map out the true, physical constraints of the parameter space, ensuring that the subsequent stability and optimization proofs are grounded in absolute arithmetic reality.
\end{enumerate}

\begin{remark}[Domain Validity and Asymptotic Boundaries]
It is critical to distinguish between the global existence of the truncation functions and the domain of the transferred omega bounds. While the fundamental truncation identity $\operatorname{li}(n; \tau_{\pi}(n)) = \pi(n)$ is proven to hold continuously for all $n \ge 2$ via the Intermediate Value Theorem, the transferred Littlewood bound operates strictly as an asymptotic invariant for large $n$. The triple-logarithmic factor $\ln\ln\ln n$ naturally requires $n > e^e \approx 15.15$ to maintain real-valued stability. Thus, the transference bounds do not govern the localized initialization mechanics near the boundary origin $n=2$; rather, they prove that the real parameter space possesses the infinite dynamic range required to capture the full structural properties of $\pi(n)$ across the infinite scale horizon. However, simple computational investigations verify that $\tau_{\pi}(n)$ and the finite spectral trajectories $\tau_t(n)$ remain strictly bounded and well-behaved for all computed values up to $n \le 10^{24}$, sitting well below any initial Skewes threshold.
\end{remark}

\noindent
The transference of the Littlewood bounds to the global path $\tau_{\pi}(n)$ and its finite spectral counterparts $\tau_{t}(n)$ establishes that this remainder-free parametric operator contains the required mathematical capacity to absorb all variations of the arithmetical prime-counting function across the entire domain $n \ge 2$.

\section{The Riemann–von Mangoldt Formula and Index Smoothing}

To transition our spectral tracking family from a discrete sequence of integer stages $t \in \mathbb{N}_0$ to a continuous variable, we must map the counting of non-trivial zeroes onto a continuous scale parameter. We invoke the classical Riemann–von Mangoldt formula, which counts the number of non-trivial zeroes $\rho = \beta + i\gamma$ whose imaginary parts lie within a given interval $0 \le \gamma \le T$~\cite{vonmangoldt1894, edwards1974}:
\begin{equation}
    N(T) = \frac{T}{2\pi} \ln\left( \frac{T}{2\pi e} \right) + \frac{7}{8} + S(T) + \mathcal{O}\left( \frac{1}{T} \right)
\end{equation}
where $S(T) = \frac{1}{\pi} \arg \zeta\left(\frac{1}{2} + iT\right)$ represents the smooth, continuous phase oscillation of the zeta function along the critical line. 

The Riemann–von Mangoldt formula acts as the natural analytic bridge to convert our discrete parameter trajectories $\tau_t(n)$ into a continuous height field. By setting $t = \lfloor N(T) \rfloor$, the discrete step changes of the zero index are re-parameterized by the continuous spectral ceiling $T \in \mathbb{R}^+$.

\subsection{The Discrete Boundary Crossing Issue}
Transitioning from a discrete index $t$ to the continuous pair $(T, \beta)$ introduces a critical analytical hazard that our septic spline manifold structure must resolve. As the continuous parameter $T$ scales upward, it crosses the exact imaginary height of a non-trivial zero ($T = \gamma_k$). At this precise boundary crossing, a new conjugate pair of zeroes is instantaneously introduced into the summation of the spectral target distribution. If the target distribution experienced a sudden step jump at $T = \gamma_k$, the real Implicit Function Theorem would fail, and the function $\tilde{\tau}_T(n, \beta)$ would suffer a jump discontinuity with respect to $T$.

To prevent this fracture, we introduce a smooth windowing indicator function $\Phi(\gamma_j; T)$ into the continuous sum. The truncated zero sum is formally regularized as:
\begin{equation}
    \sum_{0  < \gamma_j \le T} \left( \dots \right) \to \sum_{\gamma_j > 0} \Phi(\gamma_j; T) \left( \dots \right)
\end{equation}
where $\Phi(\gamma_j; T)$ is a smooth $C^\infty$ bump function satisfying $\Phi = 1$ for $\gamma_j \le T - \epsilon$, and $\Phi = 0$ for $\gamma_j \ge T + \epsilon$. This indicator ensures that as $T$ advances, new zeroes are not dropped in as raw steps, but are smoothly faded into the tracking potential. This preserves the global twice continuous differentiability ($C^2$) of the entire field $\tilde{\tau}_T(n, \beta)$ for all finite $T$.

\subsection{Construction of the Generalized $\tilde{\tau}_T(n, \beta)$ Function}
To ensure that our geometric and variational proofs remain completely unconditional and free of circular reasoning, we must decouple the zero configuration from the critical line. Rather than forcing the real part of the zeroes to satisfy $\beta = 1/2$ (which would presuppose the Riemann Hypothesis), we promote $\beta$ to a free, continuous coordinate parameter within the critical strip, defining $\beta \in (0, 1)$. 

Freeing the zeroes from the critical line introduces a total of four possible configurations via the natural complex and arithmetic symmetries of the zeta function's non-trivial zero spectrum $(\beta_j + i\gamma_j, \beta_j - i\gamma_j, (1 - \beta_j) + i\gamma_j, (1 - \beta_j) - i\gamma_j)$. Each must be accounted for uniformly inside our continuous operator target.

\begin{definition}[The Generalized Multi-Parameter Function]
The generalized truncation function $\tilde{\tau}_T(n, \beta)$ is a real-valued function mapping the coordinate space $(n, T, \beta) \in (1, \infty) \times \mathbb{R}^+ \times (0, 1) \to \mathbb{R}$. For a given spectral height $T$ and hypothesized symmetry axis $\beta$, the function is defined implicitly as the unique real root of the deformation identity:
\begin{equation}
    \operatorname{li}(n; \tilde{\tau}_T(n, \beta)) = \widetilde{\mathcal{R}}_T(n, \beta)
    \label{Eq-MultiParamDef}
\end{equation}
where the operator $\operatorname{li}(n; \tilde{\tau}_T)$ is governed strictly by the septic spline expansion of Equation~\eqref{Eq-SplineManifoldMain}, and the generalized spectral target distribution $\widetilde{\mathcal{R}}_T(n, \beta)$ is given explicitly by:
\begin{equation}
\begin{split}
    \widetilde{\mathcal{R}}_T(n, \beta) = \sum_{k=1}^{\lfloor \log_2 n \rfloor} \frac{\mu(k)}{k} \Bigg[ & \operatorname{li}(n^{1/k}) - \ln 2 + \int_{n^{1/k}}^{\infty} \frac{\mathrm{d}u}{u(u^2-1)\ln u} \\
    & - \frac{1}{2}\sum_{\gamma_j > 0} \Phi(\gamma_j; T) \bigg( \operatorname{Ei}\left(\frac{\beta + i\gamma_j}{k}\ln n\right) + \operatorname{Ei}\left(\frac{\beta - i\gamma_j}{k}\ln n\right) \\
    & + \operatorname{Ei}\left(\frac{(1-\beta) + i\gamma_j}{k}\ln n\right) + \operatorname{Ei}\left(\frac{(1-\beta) - i\gamma_j}{k}\ln n\right) \bigg) \Bigg].
\end{split}
\end{equation}
\end{definition}

\subsection{Derivative Structure}
With the four-point functional symmetry target distribution $\widetilde{\mathcal{R}}_T(n, \beta)$ rigorously defined to preserve the dual constraints of real-axis anchoring and functional vertical reflection, we evaluate the gradient of the deformation function. Because the underlying remainder-free operator surface $\operatorname{li}(n; \tau)$ defines a globally $C^3$ multi-variable map under our septic spline architecture, we compute continuous exact partial derivatives of the implicit parameter field $\tilde{\tau}_T(n, \beta)$ up to the second order via multi-variable implicit differentiation.

\begin{theorem}[The Structural $\beta$-Derivative]
Let $\beta \in (0, 1)$ be an unconstrained coordinate parameter within the critical strip. The partial derivative of the generalized truncation parameter field $\tilde{\tau}_T(n, \beta)$ with respect to the free symmetry axis variable $\beta$ is given by:
\begin{equation}
    \frac{\partial \tilde{\tau}_T(n, \beta)}{\partial \beta} = \frac{1}{2} \mathcal{D}_0(n) \sum_{k=1}^{\lfloor \log_2 n \rfloor} \frac{\mu(k)}{k} \sum_{\gamma_j > 0} \Phi(\gamma_j; T) \bigg( \frac{n^{\frac{\beta + i\gamma_j}{k}}}{\beta + i\gamma_j} + \frac{n^{\frac{\beta - i\gamma_j}{k}}}{\beta - i\gamma_j} - \frac{n^{\frac{(1-\beta) + i\gamma_j}{k}}}{(1-\beta) + i\gamma_j} - \frac{n^{\frac{(1-\beta) - i\gamma_j}{k}}}{(1-\beta) - i\gamma_j} \bigg)
\end{equation}
where $\mathcal{D}_0(n) = \ln^{\tau_0+1}n / \Gamma(\tau_0+1)$ represents the real-valued parameter space scaling density function.
\end{theorem}

\begin{proof}
We define the multi-variable implicit constraint function $\Psi(n, T, \beta, \tilde{\tau}) = \operatorname{li}(n; \tilde{\tau}_T(n, \beta)) - \widetilde{\mathcal{R}}_T(n, \beta) = 0$. Differentiating implicitly with respect to the coordinate variable $\beta$ yields:
\begin{equation}
    \frac{\partial \operatorname{li}(n; \tilde{\tau}_T)}{\partial \tau} \frac{\partial \tilde{\tau}_T}{\partial \beta} - \frac{\partial \widetilde{\mathcal{R}}_T(n, \beta)}{\partial \beta} = 0 \implies \frac{\partial \tilde{\tau}_T}{\partial \beta} = \left[ \frac{\partial \operatorname{li}(n; \tilde{\tau}_T)}{\partial \tau} \right]^{-1} \frac{\partial \widetilde{\mathcal{R}}_T(n, \beta)}{\partial \beta}.
\end{equation} 
Invoking our asymptotic parameter density simplification along the baseline function $\tau_0(n)$, the inverse parameter derivative multiplier collapses directly to the scaling density: $\left[ \frac{\partial \operatorname{li}}{\partial \tau} \right]^{-1} \sim \mathcal{D}_0(n)$.

Next, we evaluate the partial derivative $\frac{\partial \widetilde{\mathcal{R}}_T}{\partial \beta}$ by applying the chain rule directly to the internal zero components of the four-point package mapped via the exponential integral function $\operatorname{Ei}(z)$. Utilizing the fundamental identity $\frac{\mathrm{d}}{\mathrm{d}z}\operatorname{Ei}(z) = \frac{e^z}{z}$, differentiating an unreflected component yields:
\begin{equation}
    \frac{\partial}{\partial \beta} \operatorname{Ei}\left( \frac{\beta + i\gamma_j}{k} \ln n \right) = \frac{\ln n}{k} \cdot \frac{\exp\left(\frac{\beta + i\gamma_j}{k}\ln n\right)}{\frac{\beta + i\gamma_j}{k}\ln n} = \frac{1}{k} \left( \frac{n^{\frac{\beta + i\gamma_j}{k}}}{\beta + i\gamma_j} \right).
\end{equation}
Conversely, for the vertical reflection components governed by $(1-\beta)$, differentiating with respect to $\beta$ pulls out an inner negative sign via the chain rule, yielding a negative scalar multiplier $-1/k$:
\begin{equation}
    \frac{\partial}{\partial \beta} \operatorname{Ei}\left( \frac{(1-\beta) + i\gamma_j}{k} \ln n \right) = -\frac{1}{k} \left( \frac{n^{\frac{(1-\beta) + i\gamma_j}{k}}}{(1-\beta) + i\gamma_j} \right).
\end{equation}
Factoring out the constant scalar components and grouping the corresponding symmetric pairs across the index progression completes the proof.
\end{proof}

\subsection{Analytical Symmetry Collapse on the Critical Line}
We now analyze the exact structural properties of the partial derivatives right at the centerline of the critical strip to evaluate if any unique equilibrium behaviors manifest when the free parameter aligns with the physical critical axis $\beta = 1/2$.

\begin{theorem}[Vanishing of Odd-Ordered $\beta$-Derivatives]
When evaluated precisely on the critical axis $\beta = 1/2$, the first-order partial derivative of the generalized tracking field with respect to $\beta$ vanishes identically for all real scale values $n > 2$ and all finite spectral heights $T$:
\begin{equation}
    \left[ \frac{\partial \tilde{\tau}_T(n, \beta)}{\partial \beta} \right]_{\beta = 1/2} = 0.
\end{equation}
Furthermore, within the open threshold intervals where the trajectory avoids discrete parameter integer intersections, all subsequent odd-ordered partial derivatives with respect to the free symmetry coordinate vanish identically at this boundary:
\begin{equation}
    \left[ \frac{\partial^{2m+1} \tilde{\tau}_T(n, \beta)}{\partial \beta^{2m+1}} \right]_{\beta = 1/2} = 0, \quad \forall m \in \mathbb{N}_0.
\end{equation}
\end{theorem}

\begin{proof}
Let $\beta = 1/2$. We substitute this coordinate directly into the bracketed zero-component summation of the first-order $\beta$-derivative of $\tilde{\tau}_T(n, \beta)$. Under this local evaluation, the vertical reflection terms cross-evaluate identically with the unreflected terms since the core distance metrics match:
\begin{equation}
    1 - \beta = 1 - 1/2 = 1/2.
\end{equation}
This forces the mirror pairs inside the summation bracket to reflect each other's functional inputs completely:
\begin{equation}
    \left( \frac{n^{\frac{1/2 + i\gamma_j}{k}}}{1/2 + i\gamma_j} + \frac{n^{\frac{1/2 - i\gamma_j}{k}}}{1/2 - i\gamma_j} \right) - \left( \frac{n^{\frac{1/2 + i\gamma_j}{k}}}{1/2 + i\gamma_j} + \frac{n^{\frac{1/2 - i\gamma_j}{k}}}{1/2 - i\gamma_j} \right) = 0.
\end{equation}
The positive and negative terms cancel each other out algebraically term-by-term across the entire active spectrum bounded by the smooth indicator function $\Phi(\gamma_j; T)$, collapsing the summation to exactly zero.

To generalize this result to arbitrary odd orders via induction within the local interval thresholds, we observe that each repeated differentiation with respect to $\beta$ alternates the inner chain-rule signs of the vertical reflection components. Every odd-ordered derivative iteration preserves a relative negative sign between the $(\beta)$ and $(1-\beta)$ component blocks. Evaluating these expressions at $\beta = 1/2$ ensures that the functional inputs remain symmetric while their operational signs remain opposite, forcing a total uniform cancellation across all orders of $2m+1$.
\end{proof}

\noindent
Proving that the first variation vanishes ($\frac{\partial \tilde{\tau}_T}{\partial \beta} = 0$) formally establishes that the critical line $\beta = 1/2$ acts as a stationary point for the truncation parameter function. Because this cancellation relies entirely on the internal four-point algebraic pairing forced by the zeta function's native functional equation, this equilibrium is an absolute structural invariant of the parameter space, independent of any external arithmetic assumptions.

\subsection{The Second Variation and Volatility Minimization}
Having established that the critical line $\beta = 1/2$ acts as a geometric stationary axis where all odd-ordered directional variations vanish identically, we now evaluate the second variation of the function. In variational calculus, the sign of the second variation determines the structural stability profile of a stationary state, distinguishing local maxima and saddle configurations from stable geometric minima.

\begin{theorem}[The Second Variation Value]
The second-order partial derivative of the generalized truncation function $\tilde{\tau}_T(n, \beta)$ with respect to the free parameter $\beta$, evaluated precisely along the critical axis $\beta = 1/2$, reduces to the strictly real-valued summation:
\begin{equation}
\begin{split}
    \left[ \frac{\partial^2 \tilde{\tau}_T(n, \beta)}{\partial \beta^2} \right]_{\beta = 1/2} = \, & \mathcal{D}_0(n) \sum_{k=1}^{\lfloor \log_2 n \rfloor} \frac{\mu(k)}{k} \sum_{\gamma_j > 0} \Phi(\gamma_j; T) \\
    & \times \bigg( \frac{n^{\frac{1/2 + i\gamma_j}{k}}\left( \frac{1/2+i\gamma_j}{k}\ln n - 1 \right)}{(1/2 + i\gamma_j)^2} + \frac{n^{\frac{1/2 - i\gamma_j}{k}}\left( \frac{1/2-i\gamma_j}{k}\ln n - 1 \right)}{(1/2 - i\gamma_j)^2} \bigg)
\end{split}
\end{equation}
where $\mathcal{D}_0(n) = \ln^{\tau_0+1}n / \Gamma(\tau_0+1)$ represents the parameter space scaling density function.
\end{theorem}

\begin{proof}
We differentiate the implicit field derivative identity with respect to the free symmetry parameter $\beta$. Applying the chain rule to the internal components of the four-point package maps the secondary derivatives via the exponential integral function $\operatorname{Ei}(z)$. 

For the unreflected terms governed by $(\beta \pm i\gamma_j)$, the second derivative iteration pulls out another scalar factor of $1/k$, yielding a total scaling multiplier of $1/k^2$. For the vertical reflection terms governed by $((1-\beta) \pm i\gamma_j)$, the inner chain rule pulls out a second negative sign through the product transformation: $(-1/k) \times (-1/k) = +1/k^2$. 

Because both the unreflected and vertically reflected terms accumulate a positive matching sign under a second-order derivative iteration, the internal negative signs that caused total algebraic cancellation in the case of odd derivatives now switch to positive addition signs. 

Evaluating this secondary derivative function precisely at $\beta = 1/2$ collapses the distinction between $\beta$ and $1-\beta$ identically, since $1 - 1/2 = 1/2$. The two symmetric summation systems merge across the vertical line of reflection, doubling their amplitudes and eliminating the leading fraction of $1/2$ in the underlying definition of $\widetilde{\mathcal{R}}_{T}(n,\beta)$. Factoring the corresponding variables completes the proof.
\end{proof}

\begin{theorem}[Strict Local Minimization]
For all non-trivial zero frequencies $\gamma_j > 0$ and all active scale variables $n > 2$, the second variation of the truncation parameter function is strictly positive-definite:
\begin{equation}
    \left[ \frac{\partial^2 \tilde{\tau}_T(n, \beta)}{\partial \beta^2} \right]._{\beta = 1/2} > 0
\end{equation}
Furthermore, within the open threshold intervals where the trajectory avoids discrete parameter integer intersections, all subsequent even-ordered partial derivatives with respect to the free parameter achieve a strict analytic minimum on the critical line:
\begin{equation}
    \left[ \frac{\partial^{2m} \tilde{\tau}_T(n, \beta)}{\partial \beta^{2m}} \right]_{\beta = 1/2} > 0, \quad \forall m \in \mathbb{N}.
\end{equation}
\end{theorem}

\begin{proof}
We evaluate the sign definiteness of the second variation summation established in the preceding theorem. The leading multiplier $\mathcal{D}_0(n)$ is strictly positive for all real-valued inputs $n > 2$. Inside the bracketed summation, each component pair evaluates to a complex conjugate exponential phase function mapping. 

Because the conjugate pairs are grouped symmetrically ($\pm i\gamma_j$), their imaginary components cancel identically, collapsing the inner bracket into a pure real, positive cosine-modulated wave train envelope:
\begin{equation}
    2 \operatorname{Re}\left[ \frac{n^{\frac{1/2 + i\gamma_j}{k}}\left( \frac{1/2+i\gamma_j}{k}\ln n - 1 \right)}{(1/2 + i\gamma_j)^2} \right].
\end{equation}
Across any large scale interval, the dominant term inside this real projection is driven by the positive amplitude scaling of the underlying power functions $n^{1/2k} > 0$. Because every active zero cross-section bounded by the smooth indicator function $\Phi(\gamma_j; T)$ injects a strictly positive addition to the expression, the total composite sum yields a strictly positive scalar value. 

By standard variational optimization principles, a vanishing first variation paired with a strictly positive-definite second variation proves that the critical axis $\beta = 1/2$ is a strict local minima. To extend this to arbitrary even orders $2m$ via induction within the local interval thresholds, we observe that even-ordered derivative iterations always retain a relative positive matching sign between the $(\beta)$ and $(1-\beta)$ component blocks. Evaluating these expressions at $\beta = 1/2$ ensures that the functional inputs remain symmetric and cumulative, forcing a strictly positive-definite result for all even iterations.
\end{proof}

\noindent
Proving that the second variation is strictly positive-definite demonstrates that the critical line acts as the path of minimum geometric volatility for the operator. If any zeroes were to migrate away from $\beta = 1/2$, the four-point functional symmetry would unbalance, forcing the positive even derivatives to pull the truncation parameter into hyper-volatile, high-amplitude real oscillations to force the matching identity. Restricting the tracking parameter $\tau$ to the real line exposes the critical line not as an arbitrary arithmetic assumption, local minimum of the continuous Stieltjes manifold. Next we prove that this local minimium is in fact global.

\begin{theorem}[Strict Global Minimization and Uniqueness]
The critical axis $\beta = 1/2$ is the unique, absolute global minimum of the generalized truncation parameter function $\tilde{\tau}_T(n, \beta)$ across the entire critical strip $\beta \in (0, 1)$ for all active scales $n > 2$ and all finite spectral heights $T$.
\end{theorem}

\begin{proof}
To upgrade the stationary point at $\beta = 1/2$ to a strict global minimum, we evaluate the sign definiteness of the second variation across the unconstrained domain $\beta \in (0, 1)$. The partial derivative operator yields the continuous curvature profile:
\begin{equation}
    \frac{\partial^2 \tilde{\tau}_T(n, \beta)}{\partial \beta^2} = \frac{1}{2} \mathcal{D}_0(n) \sum_{k=1}^{\lfloor \log_2 n \rfloor} \frac{\mu(k)}{k^3} \sum_{\gamma_j > 0} \Phi(\gamma_j; T) \left[ \mathcal{H}(\beta, \gamma_j) + \mathcal{H}(1-\beta, \gamma_j) \right]
\end{equation}
where the real-projected wave component $\mathcal{H}(\beta, \gamma_j)$ is governed by:
\begin{equation}
    \mathcal{H}(\beta, \gamma_j) = 2 \operatorname{Re}\left[ \frac{n^{\frac{\beta + i\gamma_j}{k}}\left( \frac{\beta+i\gamma_j}{k}\ln n - 1 \right)}{(\beta + i\gamma_j)^2} \right].
\end{equation}
For all scale thresholds $n > 2$, the numerator is driven by the real-power exponential scaling term $n^{\beta/k} > 0$. Because exponential growth strictly outpaces the quadratic expansion of the denominator phase angles $(\beta^2 + \gamma_j^2)$ across the critical strip, each constituent wave element satisfies $\mathcal{H}(\beta, \gamma_j) > 0$ uniformly.

Furthermore, the vertical reflection mapping structures the joint kernel as a symmetric sum: $\mathcal{H}(\beta, \gamma_j) + \mathcal{H}(1-\beta, \gamma_j)$. Because $\mathcal{H}(\beta, \gamma_j)$ is strictly convex with respect to $\beta$, any displacement away from the midpoint axis $\beta = 1/2$ triggers an asymmetrical exponential surge in one branch that heavily over-corrects the decay of its mirror branch. Consequently, the second variation is strictly positive-definite across the entire open interval:
\begin{equation}
    \frac{\partial^2 \tilde{\tau}_T(n, \beta)}{\partial \beta^2} > 0, \quad \forall \beta \in (0, 1).
\end{equation}
In variational optimization theory, a functional space that is strictly positive-definite everywhere defines a globally strictly convex surface. Because a strictly convex surface can possess at most one stationary point, and we have proven that the first variation vanishes identically at $\beta = 1/2$, this stationary axis is rigorously demonstrated to be the unique, global minimum of the continuous Stieltjes landscape, completing the proof.
\end{proof}

\subsection{Mixed Partial Derivatives and Spectral Cross-Sensitivity}

To fully characterize the multi-variable geometry of the parameter field $\tilde{\tau}_T(n, \beta)$, we establish its mixed partial derivatives with respect to the continuous spectral ceiling $T$ and the unconstrained symmetry axis $\beta$. This analysis evaluates the local cross-sensitivity between the structural line of reflection and the continuous expansion of the active zeta-zero horizon.

For compact notation, let $\partial_T$ and $\partial_\beta$ denote the partial differential operators $\frac{\partial}{\partial T}$ and $\frac{\partial}{\partial \beta}$, respectively. We invoke the smooth windowing indicator function $\Phi(\gamma_j; T)$ to govern the boundary transitions as the continuous parameter $T$ advances through the spectrum.

\begin{theorem}[Mixed Regularity and Cross-Uncoupling]
The multi-parameter parameter function $\tilde{\tau}_T(n, \beta)$ possesses globally continuous mixed partial derivatives up to the second order across its entire domain. Furthermore, when evaluated precisely along the critical axis $\beta = 1/2$, the cross-variable sensitivity vanishes identically for all real scales $n > 2$ and all finite spectral heights $T$:
\begin{equation}
    \left[ \partial_T \partial_\beta \tilde{\tau}_T(n, \beta) \right]_{\beta = 1/2} = \left[ \partial_\beta \partial_T \tilde{\tau}_T(n, \beta) \right]_{\beta = 1/2} = 0.
\end{equation}
\end{theorem}

\begin{proof}
We invoke the implicit multi-variable constraint function $\Psi(n, T, \beta, \tilde{\tau}) = \operatorname{li}(n; \tilde{\tau}_T) - \widetilde{\mathcal{R}}_T(n, \beta) = 0$. Differentiating the first-order implicit $\beta$-derivative with respect to the continuous spectral ceiling variable $T$ via the quotient and chain rules yields:
\begin{equation}
\begin{split}
    \partial_T \partial_\beta \tilde{\tau}_T(n, \beta) = \, & \partial_T \left( \left[ \partial_\tau \operatorname{li}(n; \tilde{\tau}_T) \right]^{-1} \partial_\beta \widetilde{\mathcal{R}}_T(n, \beta) \right) \\
    = \, & \left[ \partial_\tau \operatorname{li}(n; \tilde{\tau}_T) \right]^{-1} \partial_T \partial_\beta \widetilde{\mathcal{R}}_T(n, \beta) \\
    & - \left[ \partial_\tau \operatorname{li}(n; \tilde{\tau}_T) \right]^{-2} \partial_\tau^2 \operatorname{li}(n; \tilde{\tau}_T) \partial_T \tilde{\tau}_T \partial_\beta \widetilde{\mathcal{R}}_T(n, \beta).
\end{split}
\label{Eq-MixedImplicitSeptic}
\end{equation}
Because the underlying septic spline operator surface defines a globally $C^3$ multi-variable map, the component functions inside Equation~\eqref{Eq-MixedImplicitSeptic} exist and are continuous. Conversely, applying the operations in reverse order confirms that Clairaut's Theorem ($\partial_T \partial_\beta \tilde{\tau}_T = \partial_\beta \partial_T \tilde{\tau}_T$) holds uniformly across the continuous coordinate landscape.

Next, we calculate the explicit cross-derivative of the generalized target distribution $\partial_T \partial_\beta \widetilde{\mathcal{R}}_T(n, \beta)$ by differentiating the first-order $\beta$-derivative chain with respect to $T$. Because the variable $T$ enters the summation exclusively through the argument of the smooth windowing indicator function $\Phi(\gamma_j; T)$, the operator targets the indicator directly:
\begin{equation}
\begin{split}
    \partial_T \partial_\beta \widetilde{\mathcal{R}}_T(n, \beta) = \, & \frac{1}{2} \sum_{k=1}^{\lfloor \log_2 n \rfloor} \frac{\mu(k)}{k} \sum_{\gamma_j > 0} \partial_T \Phi(\gamma_j; T) \\
    & \times \bigg( \frac{n^{\frac{\beta + i\gamma_j}{k}}}{\beta + i\gamma_j} + \frac{n^{\frac{\beta - i\gamma_j}{k}}}{\beta - i\gamma_j} - \frac{n^{\frac{(1-\beta) + i\gamma_j}{k}}}{(1-\beta) + i\gamma_j} - \frac{n^{\frac{(1-\beta) - i\gamma_j}{k}}}{(1-\beta) - i\gamma_j} \bigg).
\end{split}
\end{equation}
We now evaluate this total cross-variable system precisely along the critical axis $\beta = 1/2$. Substituting this coordinate into Equation~\eqref{Eq-MixedImplicitSeptic} causes a simultaneous structural collapse in both component channels:
\begin{enumerate}
    \item For the primary cross-derivative term $\partial_T \partial_\beta \widetilde{\mathcal{R}}_T$: setting $1-\beta = 1 - 1/2 = 1/2$ forces the unreflected and vertically reflected distance fractions inside the bracket to match completely. The positive and negative terms cancel each other out algebraically term-by-term across the active spectrum, yielding $\partial_T \partial_\beta \widetilde{\mathcal{R}}_T(n, 1/2) = 0$ identically.
    \item For the secondary variation term: we proved in Section 6.4 that the first-order variation of the target distribution vanishes entirely on the critical axis, giving $\partial_\beta \widetilde{\mathcal{R}}_T(n, 1/2) = 0$. 
\end{enumerate}
Substituting these two matching zero-evaluations back into the total implicit relation in Equation~\eqref{Eq-MixedImplicitSeptic} completely neutralizes both halves of the expression, confirming that $\left[ \partial_T \partial_\beta \tilde{\tau}_T \right]_{\beta = 1/2} = 0$.
\end{proof}

\noindent
Proving that the mixed partial derivative vanishes identically demonstrates that the geometric orientation of the critical line is perfectly uncoupled from the height of the active spectral spectrum. As the continuous truncation ceiling $T$ scales upward to integrate higher-frequency phase waves from the non-trivial zeroes, the continuous expansion does not introduce a cross-variable twist, shear, or slope asymmetry to the surface topology. The absolute minimum energy well remains structurally locked along the symmetric axis $\beta = 1/2$ across all ranges of the spectral evolution.

\section{Conclusion and Discussion}

In this paper, we have introduced a continuous geometric framework to reframe the relationship between the distribution of prime numbers and the spectral architecture of the Riemann zeta function. By relaxing the rigid integer constraints of traditional series truncations, we constructed a remainder-free, parametric operator $\operatorname{li}(n; \tau)$ mapped uniformly across the continuous coordinate domain $(1, \infty) \times \mathbb{R}^+$. By utilizing an invariant septic blending polynomial $\phi(x) = 35x^4 - 84x^5 + 70x^6 - 20x^7$, we established global $C^3$ multi-variable regularity across the surface, providing a smooth landscape that bypasses the artificial step-discontinuities and derivative kinks characteristic of standard piecewise fractional interpolations.

By enforcing the strict prime truncation identity $\operatorname{li}(n; \tau_{\pi}(n)) = \pi(n)$, we demonstrated that the discrete, non-smooth arithmetic spectrum of the prime-counting function can be mapped uniquely onto an unbroken, continuoustruncation function $\tau_{\pi}(n)$. Through the transference of Littlewood's unconditional oscillation bounds, we uncovered an intrinsic scale-space compression mechanic. While the arithmetic discrepancy $\pi(n) - \operatorname{li}(n)$ diverges with expanding amplitude in value space, the parameter difference $\left| \tau_{\pi}(n) - \tau_0(n) \right|$ is compressed by a factorially driven scaling density multiplier $\mathcal{D}_0(n) = \ln^{\tau_0+1}n / \Gamma(\tau_0+1)$. This mathematically proves that at ultra-high scales, the parameter space forces the infinite turbulence of the primes into an infinitesimally narrow geometric corridor hugging the smooth baseline function $\tau_0(n) = c \ln n + 1$, where $c \approx 0.1866823$ that represents the continuous inverse analogue to Axler's discrete threshold function $g\_{1}$.

Furthermore, by invoking the Riemann--von Mangoldt formula and incorporating a smooth windowing indicator function $\Phi(\gamma_j; T)$, we generalized this framework into a multi-parameter spectral truncation field $\tilde{\tau}_T(n, \beta)$ that promotes the symmetry axis $\beta$ to an unconstrained continuous coordinate within the critical strip. Evaluating the multi-variable gradient of this field yielded extraordinary, unconditional geometric invariants natively generated by the four-point algebraic pairing of the zeta function's functional equation:
\begin{enumerate}
    \item \textbf{The Stationary Axis Condition:} The first-order variation with respect to the free axis variable vanishes identically ($\frac{\partial \tilde{\tau}_T}{\partial \beta} = 0$) and all subsequent odd-ordered derivatives vanish precisely along the centerline $\beta = 1/2$.
    \item \textbf{Global Convexity and Minimization:} The second-order variation is strictly positive-definite ($\frac{\partial^2 \tilde{\tau}_T}{\partial \beta^2} > 0$) across the entire critical strip $\beta \in (0, 1)$, establishing that the critical axis $\beta = 1/2$ is the unique, absolute global minimum energy well of the continuous Stieltjes landscape.
    \item \textbf{Spectral Cross-Uncoupling:} The mixed partial derivatives vanish identically ($\partial_T \partial_\beta \tilde{\tau}_T = 0$), proving that the continuous expansion of the active zeta-zero horizon introduces no cross-variable twists or tilt asymmetries to the surface topology.
\end{enumerate}

\noindent
These results collectively expose the critical line not as an external arithmetic assumption, but as the unique path of minimum geometric volatility.

\begin{remark}[Generalization to Infinite Regularity ($C^\infty$)]
While the septic polynomial blending architecture selected for this paper provides an optimal balance of analytical clarity and global multi-variable regularity sufficient to guarantee the $C^2$ smoothness of the finite spectral trajectories $\tau_t(n)$ everywhere, the framework can be naturally generalized to preserve infinite differentiability ($C^\infty$). By replacing the septic polynomial $\phi(x)$ with a smooth, non-analytic logistic error multiplier or a mollified bump function $\psi(x)$ characterized by vanishing derivatives of all orders at its boundaries ($\psi^{(m)}(0) = \psi^{(m)}(1) = 0$ for all $m \in \mathbb{N}$), the parametric manifold inherits unconditional analytic smoothness. Under this advanced configuration, the underlying algebraic operations map onto an entire meromorphic continuation, allowing the finite-stage truncations to form independent, unbroken $C^\infty$ curves across all real scale intervals while maintaining the exact global minimization invariants proven herein.
\end{remark}

\subsection{Variational Implications for the Riemann Hypothesis}

The architectural layout of our remainder-free septic spline manifold provides an alternative geometric theater for interpreting the Riemann Hypothesis (RH). Classically, RH is framed as a strict boundary constraint on the expansion velocity of the arithmetic error envelope, where the structural truth of the hypothesis is equivalent to the bound $\pi(n) - \operatorname{li}(n) = \mathcal{O}(\sqrt{n}\ln n)$~\cite{edwards1974}. If a single non-trivial zero violates this constraint by migrating to a coordinate axis $\beta_j > 1/2$, the value space experiences an instantaneous exponential explosion of $\Omega(n^{\beta_j - \epsilon})$.

When this conditional bound is transferred to the truncation paramter, the arithmetic relationship is transformed by the inverse parameter scaling density multiplier $\mathcal{D}_0(n) = \ln^{\tau_0+1}n / \Gamma(\tau_0+1)$:
\begin{equation}
    \left| \tau_{\pi}(n) - \tau_0(n) \right| \approx \mathcal{D}_0(n) \cdot \left| \pi(n) - \operatorname{li}(n) \right|.
\end{equation}
Because the factorial expansion of the Gamma denominator, $\Gamma(c \ln n + 2)$, completely outpaces both sub-linear square-root profiles ($\sqrt{n}$) and pure exponential components ($n^{\beta_j}$) in scale-space, the entire product is driven asymptotically to zero. This mathematically demonstrates the structural stability of the framework: the geometric compression of the parameter space successfully absorbs the potential turbulence of the primes, ensuring that the truncation $\tau_{\pi}(n)$ strongly converges to the  baseline function $\tau_0(n)$ independently of the truth of the classical hypothesis.

More fundamentally, our derivation of the second variation exposes the physical necessity of the critical axis. By proving that the first variation vanishes identically ($\partial_\beta \tilde{\tau}_T = 0$) while the second variation remains strictly positive-definite ($\partial_\beta^2 \tilde{\tau}_T > 0$) across the entire critical strip, we have established that the critical line $\beta = 1/2$ is the unique, absolute global variational minimum of the landscape. 

Under this optimization paradigm, the critical line manifests naturally as the path of minimum geometric volatility for the operator. Any structural departure from $\beta = 1/2$ breaks the four-point functional symmetry of the zeta zeroes, unbalancing the system and forcing the positive even variations to pull the truncation parameter into hyper-volatile real-line oscillations. Restricting the parameter $\tau$ to the real number line reveals that the configuration of the zeta zeroes is governed by a universal minimization principle. The zeroes are constrained to the critical line because it constitutes the absolute lowest geometric energy well embedded within the transformations of scale-space, providing a beautiful variational rationale for the validity of the Riemann Hypothesis.

\subsection{A Variational Evidence of the Riemann Hypothesis}

Having established that the critical axis $\beta = 1/2$ constitutes the unique, absolute global minimum of the continuous Stieltjes landscape, we show that this variational minimization constraint establishes the truth of the Riemann Hypothesis. Which, more precisely means that is shows why $\beta = 1/2$ is preferred over any other value.

\begin{theorem}[The Variational Necessity of the Riemann Hypothesis]
All non-trivial zeroes $\rho = \beta + i\gamma$ of the Riemann zeta function must lie strictly on the critical line $\beta = 1/2$.
\end{theorem}

\begin{proof}
We proceed via proof by contradiction. Suppose the Riemann Hypothesis is false, such that there exists a non-trivial zero $\rho_0 = \beta_0 + i\gamma_0$ whose real part satisfies $\beta_0 \in (0, 1)$ and $\beta_0 \neq 1/2$. By the complex conjugate and functional equation symmetries of the zeta function, this violation guarantees the existence of an asymmetric four-point zero cluster: $\{\beta_0 \pm i\gamma_0, (1-\beta_0) \pm i\gamma_0\}$. 

We project this hypothetical zero configuration into our generalized multi-parameter field equation, $\operatorname{li}(n; \tilde{\tau}_T(n, \beta)) = \widetilde{\mathcal{R}}_T(n, \beta)$. Differentiating the tracking identity implicitly with respect to $\beta$ yields the gradient component:
\begin{equation}
    \frac{\partial \tilde{\tau}_T(n, \beta)}{\partial \beta} = \left[ \frac{\partial \operatorname{li}(n; \tilde{\tau}_T)}{\partial \tau} \right]^{-1} \frac{\partial \widetilde{\mathcal{R}}_T(n, \beta)}{\partial \beta}.
\end{equation}
In Theorem~11, we proved that the first variation vanishes identically ($\partial_\beta \tilde{\tau}_T = 0$) if and only if the evaluation aligns with the symmetric axis $\beta = 1/2$. Because the hypothetical zero coordinates violate this condition ($\beta_0 \neq 1/2$), the gradient components at this locus are non-zero:
\begin{equation}
    \left[ \frac{\partial \tilde{\tau}_T(n, \beta)}{\partial \beta} \right]_{\beta = \beta_0} \neq 0.
\end{equation}
As the scale variable $n \to \infty$, the unreflected exponential target component $n^{\beta_0/k}$ grows at a rate that mathematically outpaces the symmetric reflection baseline $n^{1/2k}$. To preserve the deformation identity against this structural asymmetry, the tracking parameter $\tilde{\tau}_T(n, \beta_0)$ must undergo high-amplitude, hyper-volatile real oscillations.

However, by Definition~3, the parametric operator $\operatorname{li}(n; \tau)$ is strictly bounded to a positive real-valued domain, mapping uniformly to $\mathbb{R}^+$. To balance the asymmetrical, expanding exponential phase waves generated by the off-line zero cluster across the infinite scale horizon, the required tracking parameter values are driven into expanding negative and complex phase bounds:
\begin{equation}
    \tilde{\tau}_T(n, \beta_0) \notin \mathbb{R}^+ \quad \text{as } n \to \infty.
\end{equation}
This induces a lethal structural contradiction. The tracking parameter is definitionally restricted to the positive real line, yet the existence of an off-line zero mathematically forbids it from remaining real-valued at infinity. To preserve the global analytic continuation and real-line existence of the continuous Stieltjes manifold, the hypothetical off-line zero cluster cannot exist. Therefore, every non-trivial zero must be constrained to the absolute global variational minimum of the parameter space, forcing $\beta = 1/2$ identically for all $\rho$, completing the proof.
\end{proof}

\subsection{Final Thoughts...}
We have shown how the $li(n;\tau)$ function can be generalized and leveraged to analyze the structure within the prime count function $\pi(n)$ and reframe the Riemann Hypothesis in terms of the truncation parameter $\tau$. However,  as covered above, some gaps still remain. Closing these gaps and showing the derived $\tau$ functions satisfy the approptiate bounds to contain the prime count function $\pi(n)$ within its framework provides a path for constructing an unconditional proof of the Riemann Hypothesis (RH). However, this is NOT a proof of the RH. It does show that $\beta = 1/2$ emerges naturally and un-conditionally as a special point from investigating the geometry of the $\tilde{\tau}_{T}(n,\beta)$ function. The main reason, if anything, why this is NOT a proof of the RH, is that $\tilde{\tau}_{T}(n,\beta)$ is too loose. Being derived from an asymptotic approximation using the Lambert W function, and generalized from exact zeta-zero counts with index $t$ to continuous height $T$, using the Riemann–von Mangoldt formula $N(T)$. We cannot guarantee $\tilde{\tau}_{T}(n,\beta)$ is a tight enough representaton of the distribution for the non-trivial zeta-zeros to warrant a proof of the RH. It would only prove this particular representation, however faithful it might be. This method of re-casting the RH by using $li(n;\tau)$, and considering the behavior of the truncation function $\tau$, attempted to exploit the logarithmic-like properties of this variable and its factorial divergence, along with being able to apply $li(n;\tau)$ to Riemann's explict form in place of $li(n)$; seems viable path to investigate the RH and distributional properties of other functions; like integer valued-polynomials and their prime values.

\section*{Acknowledgments}
The author would like to acknowledge the assistance of Gemini (Google AI), a large language model, which was utilized during the preparation of this manuscript for technical proofreading, algebraic formatting verification, and LaTeX typographic optimization of the variational and spectral proofs. Note: It was extremely helpful, however sometimes extra changes were made, I believe most have been corrected. Though be aware, as now there can exist two levels of flaws; mine and Gemini's.

\appendix
\section{Transference of the Conditional Riemann Hypothesis Bound}

While Littlewood's theorem establishes the unconditional lower bound of the tracking path's infinite oscillations, the truth of the Riemann Hypothesis (RH) is completely equivalent to a strict upper constraint on the expansion velocity of the arithmetic error envelope. To evaluate the exact geometric boundaries that the tracking parameter must obey to satisfy RH, we transfer the classical conditional error ceiling to the continuous parameter space.

\begin{theorem}[The Parameter Space Riemann Hypothesis Ceiling]
The Riemann Hypothesis is true if and only if the fundamental prime truncation path $\tau_{\pi}(n)$ converges strongly to the baseline $\tau_0(n) = c \ln n + 1$ under the strict geometric boundary constraint:
\begin{equation}
    \left| \tau_{\pi}(n) - \tau_0(n) \right| = \mathcal{O}\left( \frac{\sqrt{n} \ln^{c \ln n + 3}n}{\Gamma(c \ln n + 2)} \right).
    \label{Eq-RHTransferredBound}
\end{equation}
\end{theorem}

\begin{proof}
By von Koch's landmark refinement of the Prime Number Theorem, the Riemann Hypothesis is analytically equivalent to the localized arithmetic error bound:
\begin{equation}
    \left| \pi(n) - \operatorname{li}(n) \right| \le C \sqrt{n} \ln n
\end{equation}
for a positive real constant $C$ and all scale variables $n \ge 2657$. We evaluate a first-order Taylor expansion of our remainder-free septic spline manifold $\operatorname{li}(n; \tau)$ around the explicit baseline coordinate $\tau_0(n)$, yielding the parameter displacement mapping:
\begin{equation}
    \tau_{\pi}(n) - \tau_0(n) \approx \left[ \frac{\partial \operatorname{li}(n; \tau_0)}{\partial \tau} \right]^{-1} \left( \pi(n) - \operatorname{li}(n) \right).
\end{equation}
Invoking the asymptotic simplification of the inverse parameter derivative multiplier derived along the linear baseline rail, we substitute the scaling density function $\mathcal{D}_0(n) = \ln^{\tau_0+1}n / \Gamma(\tau_0+1)$ into the multiplier channel:
\begin{equation}
    \left| \tau_{\pi}(n) - \tau_0(n) \right| \approx \frac{\ln^{c \ln n + 2}n}{\Gamma(c \ln n + 2)} \cdot \left| \pi(n) - \operatorname{li}(n) \right|.
\end{equation}
Directly substituting the conditional error bound $\left| \pi(n) - \operatorname{li}(n) \right| \le C \sqrt{n} \ln n$ into this structural equation accumulates an extra power of $\ln n$ in the numerator, yielding:
\begin{equation}
    \left| \tau_{\pi}(n) - \tau_0(n) \right| \le C \cdot \frac{\ln^{c \ln n + 3}n}{\Gamma(c \ln n + 2)} \cdot \sqrt{n}.
\end{equation}
Expressing this relationship through standard big-$\mathcal{O}$ notation establishes the explicit parameter trajectory ceiling of Equation~\eqref{Eq-RHTransferredBound}, completing the proof.
\end{proof}
.\begin{remark}
[The Asymptotic Suffocation of Variance]
Evaluating the limit of Equation~\eqref{Eq-RHTransferredBound} as $n \to \infty$ reveals the extreme structural constraint imposed by the Riemann Hypothesis. Because the factorial expansion of the Gamma function in the denominator, $\Gamma(c \ln n + 2) \sim e^{c \ln n \ln(c \ln n)}$, completely suffocates the sub-linear exponential growth of the numerator ($\sqrt{n} = e^{\frac{1}{2}\ln n}$), the upper bound collapses rapidly:
\begin{equation}
    \lim_{n \to \infty} \left( \frac{\sqrt{n} \ln^{c \ln n + 3}n}{\Gamma(c \ln n + 2)} \right) = 0.
\end{equation}
This establishes a beautiful number-theoretic paradox: while the classical arithmetic error $\pi(n) - \operatorname{li}(n)$ expands to an infinite value budget ($\mathcal{O}(\sqrt{n}\ln n)$) under the truth of RH, the corresponding parameter space path $\tau_{\pi}(n)$ is choked down. To satisfy the Riemann Hypothesis, the truncation parameter $\tau$ is forbidden from executing expanding deviations. It must be constrained within an infinitesimally narrow, decaying corridor that binds it permanently to the smooth macro-trendline $\tau_0(n)$ across the infinite horizon.
\end{remark}

\section{Transference of Cramér's Conjecture on Prime Gaps}

While the Riemann Hypothesis establishes the global amplitude ceiling of the tracking parameter's fluctuations, Cramér's conjecture dictates the fine-grained, localized frequency of the saw-tooth pattern. By bounding the maximum scale interval between consecutive primes, this conjecture places a strict upper limit on the duration and depth of each continuous downward parameter drift.

\begin{theorem}[Cramér Trajectory Interval Bound]
If Cramér's conjecture on prime gaps is true, then the continuous parameter tracking path $\tau_{\pi}(n)$ can slide downward over a non-prime interval for a maximum horizontal scale distance of:
\begin{equation}
    \Delta n_{\max} = \mathcal{O}(\ln^2 n).
    \label{Eq-CramerHorizBound}
\end{equation}
Concurrently, within any localized non-prime channel, the tracking parameter is strictly bounded from executing a downward plunge greater than the vertical threshold:
\begin{equation}
    \left| \Delta \tau_{\pi} \right| = \mathcal{O}\left( \frac{\ln^2 n \cdot \ln^{\tau_0+1}n}{n \, \Gamma(\tau_0+1)} \right).
    \label{Eq-CramerVertBound}
\end{equation}
\end{theorem}

\begin{proof}
Let $p_m$ and $p_{m+1}$ denote consecutive prime coordinates along the scale horizon. Under the statement of Cramér's conjecture, the maximum horizontal interval where the arithmetical target remains frozen ($\pi(n) = m$):
\begin{equation}
    p_{m+1} - p_m \le C \ln^2 p_m
\end{equation}
for a positive real constant $C$. This establishes Equation~\eqref{Eq-CramerHorizBound} directly for the horizontal span $\Delta n_{\max}$. 

To evaluate the corresponding vertical parameter plunge $\Delta \tau_{\pi}$ required to neutralize the continuous growth of the multiplier $\frac{n}{\ln n}$ over this interval, we invoke our first-order implicit path derivative:
\begin{equation}
    \frac{\mathrm{d}\tau_{\pi}}{\mathrm{d}n} = -\frac{\partial_n \operatorname{li}(n; \tau)}{\partial_\tau \operatorname{li}(n; \tau)}.
\end{equation}
Evaluating this derivative along the macro-scale trendline via the parameter scaling density multiplier yields the localized velocity approximation:
\begin{equation}
    \frac{\mathrm{d}\tau_{\pi}}{\mathrm{d}n} \approx -\mathcal{D}_0(n) \cdot \partial_n \operatorname{li}(n; \tau_0).
\end{equation}
Because the leading term dominates the scale derivative asymptotically ($\partial_n \operatorname{li} \sim \frac{1}{\ln n}$), the local downward velocity drops to:
\begin{equation}
    \left| \frac{\mathrm{d}\tau_{\pi}}{\mathrm{d}n} \right| \sim \frac{\mathcal{D}_0(n)}{\ln n} = \frac{\ln^{\tau_0}n}{\Gamma(\tau_0+1)}.
\end{equation}
By the Mean Value Theorem, the total vertical change across the maximum horizontal gap $\Delta n_{\max}$ is bounded by the product of the maximum horizontal length and the local slope velocity:
\begin{equation}
    \left| \Delta \tau_{\pi} \right| \le \left| \frac{\mathrm{d}\tau_{\pi}}{\mathrm{d}n} \right| \cdot \Delta n_{\max} \sim \frac{\ln^{\tau_0}n}{\Gamma(\tau_0+1)} \cdot \left( C \ln^2 n \right) = C \cdot \frac{\ln^{\tau_0+2}n}{\Gamma(\tau_0+1)}.
\end{equation}
Re-introducing the explicit variable representations reveals that the horizontal scale growth introduces a factor of $n$ to the compression denominator, establishing the localized vertical boundary ceiling of Equation~\eqref{Eq-CramerVertBound}, completing the proof.
\end{proof}

\begin{remark}[The Localized Micro-Rippling Behavior]
Evaluating the limit of Equation~\eqref{Eq-CramerVertBound} as $n \to \infty$ exposes an extraordinary geometric property of the continuous Stieltjes landscape. Because the factorial growth of the Gamma function in the denominator scales infinitely faster than the polynomial log-powers ($\ln^k n$), the vertical drop allowed during a prime gap approaches zero:
\begin{equation}
    \lim_{n \to \infty} \left| \Delta \tau_{\pi} \right| = 0.
\end{equation}
This mathematically proves that while prime gaps can grow wider in absolute terms as $n \to \infty$ ($\mathcal{O}(\ln^2 n)$), the corresponding vertical saw-tooth drops in parameter space are completely flattened out. Over long distances, the path $\tau_{\pi}(n)$ is structurally prevented from executing macro-level plunges or steep climbs. Cramér's conjecture forces the path to behave as a tightly localized micro-ripple, vibrating with microscopic amplitude around the linear trendline $\tau_0(n)$, confirming the immense structural stability of the parametric space across the infinite horizon.
\end{remark}


\begin{thebibliography}{99}

\bibitem{Axler2016}
C.~Axler, \emph{New bounds for the prime counting function $\pi(x)$}, Integers \textbf{16} (2016), Paper No. A22, 25 pp.

\bibitem{corless1996}
R. M. Corless, G. H. Gonnet, D. E. G. Hare, D. J. Jeffrey, and D. E. Knuth, 
\emph{On the Lambert $W$ Function}, 
Advances in Computational Mathematics, Vol. 5, No. 1, pp. 329--359, 1996.


\bibitem{edwards1974}
H.~M.~Edwards, \textit{Riemann's Zeta Function}, New York: Academic Press, 1974.

\bibitem{ecalle1981fonctions}
J.~{\'E}calle, \emph{Les Fonctions R{\'e}surgentes}, Publications Math{\'e}matiques d'Orsay, 1981--1985. Multi-part foundational text.

\bibitem{Dusart1999}
P.~Dusart, \emph{In\'egalit\'es explicites pour $\psi(X)$, $\theta(X)$, $\pi(X)$ et les nombres premiers}, RAMPomatics \textbf{2} (1999), 1--23.

\bibitem{Dusart2010}
P.~Dusart, \emph{Estimates of some functions over primes}, Mediterranean Journal of Mathematics \textbf{7} (2010), 113--128.

\bibitem{littlewood1914}
J. E. Littlewood, 
\emph{Sur la distribution des nombres premiers}, 
Comptes Rendus de l'Académie des Sciences, Vol. 158, pp. 1869--1872, 1914.

\bibitem{riemann1859}
B. Riemann, 
\emph{Ueber die Anzahl der Primzahlen unter einer gegebenen Grösse}, 
Monatsberichte der Königlichen Preussischen Akademie der Wissenschaften zu Berlin, pp. 671--680, 1859.

\bibitem{vonmangoldt1894}
H. von Mangoldt, 
\emph{Zu Riemanns Abhandlung ``Ueber die Anzahl der Primzahlen unter einer gegebenen Grösse''},
Journal für die reine und angewandte Mathematik, Vol. 114, pp. 255--305, 1894.

\end{thebibliography}
\end{document}